\documentclass[10pt,twocolumn]{article}

\usepackage{shock}
\usepackage{graphicx}
\usepackage{latexsym}
\usepackage{layout}
\usepackage{amssymb}
\usepackage{amsmath}

\input amssym.def
\input amssym.tex

\newcounter{th}
\newenvironment{theorem}{\par
\refstepcounter{th} {\sc Theorem
\arabic{th}.}\it}

\newcounter{lm}
\newenvironment{lemma}{\par \medskip \refstepcounter{lm} {\sc Lemma
\arabic{lm}.}\it}{\medskip}

\newcounter{co}

\newcounter{df}

\newcounter{nt}
\newenvironment{remark}{\par
\refstepcounter{nt} {\sc
Remark.}}{\medskip}

\newcommand{\R}{\mathbb{R}}

\newcommand{\Cb}{\mathbf{C}}
\newcommand{\C}{C}

\newcommand{\E}{\mathcal{E}}
\newcommand{\I}{\mathcal{I}}
\newcommand{\J}{\mathcal{J}}
\newcommand{\Q}{\mathcal{O}}

\newcommand{\cnst}{\mathrm{const}}
\newcommand{\rnk}{\mathop{\mathrm{rank}}\nolimits}

\newcommand{\pot}{\varphi}
\newcommand{\potv}{\psi}
\newcommand{\vel}{\mathbf{v}}
\newcommand{\nbl}{\mathop{\nabla}\nolimits}
\newcommand{\lpl}{\mathop{\Delta}\nolimits}
\newcommand{\visc}{\nu}
\newcommand{\act}{W}
\newcommand{\sact}{w}
\newcommand{\Lag}{L}
\newcommand{\Ham}{H}
\newcommand{\ham}{h}

\newcommand{\shk}{S}

\newcommand{\xb}{\mathbf{x}}
\newcommand{\ab}{\mathbf{a}}
\newcommand{\bb}{\mathbf{b}}
\newcommand{\ub}{\mathbf{u}}
\newcommand{\ui}{u}
\newcommand{\vb}{\mathbf{z}}
\newcommand{\vi}{z}
\newcommand{\pb}{\mathbf{p}}

\newcommand{\y}{y}

\newcommand{\g}{\mathbf{g}}
\newcommand{\G}{\mathbf{G}}
\newcommand{\id}{\mathrm{id}}

\newcommand{\nm}{m}
\newcommand{\na}{d}
\newcommand{\nk}{j}
\newcommand{\n}{l}
\newcommand{\nl}{n}

\newcommand{\Tr}{\Theta}
\newcommand{\Qd}{\Xi}
\newcommand{\II}{\mathrm{II}}

\newcommand{\X}{\mathcal{X}}
\newcommand{\K}{\mathcal{K}}
\newcommand{\W}{\mathcal{W}}
\newcommand{\T}{\mathcal{T}}
\newcommand{\Ab}{\mathbf{A}}
\newcommand{\Ai}{\mathcal{A}}
\newcommand{\Bb}{\mathbf{B}}

\newcommand{\df}{\mathrm{d}}

\title{Perestroikas of Shocks and Singularities of Minimum Functions}

\author{
{\sc Ilya A. Bogaevsky}\\ \\
{\small Independent University of Moscow}\\
{\small Bolsho{\u\i} Vlas$'$evski{\u\i} per.\,11, Moscow 121002, Russia}\\
{\small \it E-mail: bogaevsk@mccme.ru}}

\date{}

\begin{document}

\maketitle

\begin{abstract}\noindent
The shock discontinuities, generically present in inviscid
solutions of the forced Burgers equation, and their bifurcations
happening in the course of time (\emph{perestroikas}) are
classified in two and three dimensions -- the one-dimensional case
is well known. This classification is a result of selecting among
all the perestroikas occurring for minimum functions depending
generically on time, the ones permitted by the convexity of the
Hamiltonian of the Burgers dynamics. Topological restrictions on
the admissible perestroikas of shocks are obtained. The resulting
classification can be extended to the so-called viscosity
solutions of a Hamilton--Jacobi equation, provided the Hamiltonian
is convex.


\end{abstract}

\renewcommand{\thefootnote}{\fnsymbol{footnote}}
\footnotetext{Partially supported by INTAS-00-0259,
NWO-RFBR-047.008.005, and RFBR-02-01-00655.}

\section{Introduction}

The subject of this paper is singularities and {\it perestroikas}
($=$ bifurcations, metamorphoses, or transitions) of shocks in
plane and space (two and three diemsnions).
{\it Shocks} are discontinuities of limit
potential solutions of the Burgers equation with vanishing
viscosity and external potential force. The Burgers equation is
just the Navier--Stokes equation without the pressure term -- its
theory is well described in the survey \cite{FB00}.

If we remove the viscosity term as well, we will get the equation
describing irrotational flow of a medium consisting of
non-interacting particles moving in an external potential force
field. This is a model of formation of the large-structure of the
Universe or the so-called ``Zel$'$dovich approximation''
\cite{Zel70}, \cite{AShZ}. According to this model, after some
time the density becomes infinite because the fastest particles
overrun the slowest ones. The points with infinite density form
{\it caustics} in space which evolve with time. Generic
singularities and perestroikas of caustics are described in
\cite{ArnCat}, \cite{AGVI}.

Limit potential solutions of the Burgers equation describe motion
when the particles cannot pass through each other -- they adhere
in some sense. The adhesion is a result of interaction between
them described by the vanishing viscosity. The shocks of limit
solutions, like the above caustics, are formed by the points with
infinite density. But there is a big difference: a particle cannot
leave a shock but can leave a caustic. In other words, the matter
accumulates in the shocks. This is the so-called adhesion model of
matter evolution in the Universe describing the formation of
cellular structure of the matter described, for example, in
\cite{GMS91}.

Shocks appear even if the initial condition is {\it smooth} --
here and further this term means infinite differentiability. We
describe local singularities and perestroikas of shocks in
physically interesting dimensions $d \le 3$ with generic smooth
initial conditions. It means that we are interested only in
singularities and perestroikas which are stable with respect to
any sufficiently small perturbations of the smooth initial
condition. Other singularities and perestroikas can be killed by
its arbitrarily small perturbation.

It turns out that at typical times such a generic shock has
singularities from a finite list and at separate times experiences
perestroikas. In order to describe them we consider the {\it world
shock} lying in the space-time and formed by the points where our
limit solution is not smooth -- our shocks in space or {\it
instant shocks} are just sections of the world shock with
isochrones $t=\cnst$. The number of topologically different
generic perestroikas proves to be finite as well.

\begin{figure}[h]
\begin{center}
\includegraphics{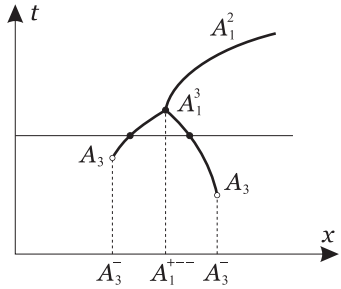}
\caption{Singularities of world shocks in plane generate
perestroikas of instant shocks on line} \label{f1a}
\end{center}
\end{figure}

\begin{figure*}
\begin{center}
\includegraphics{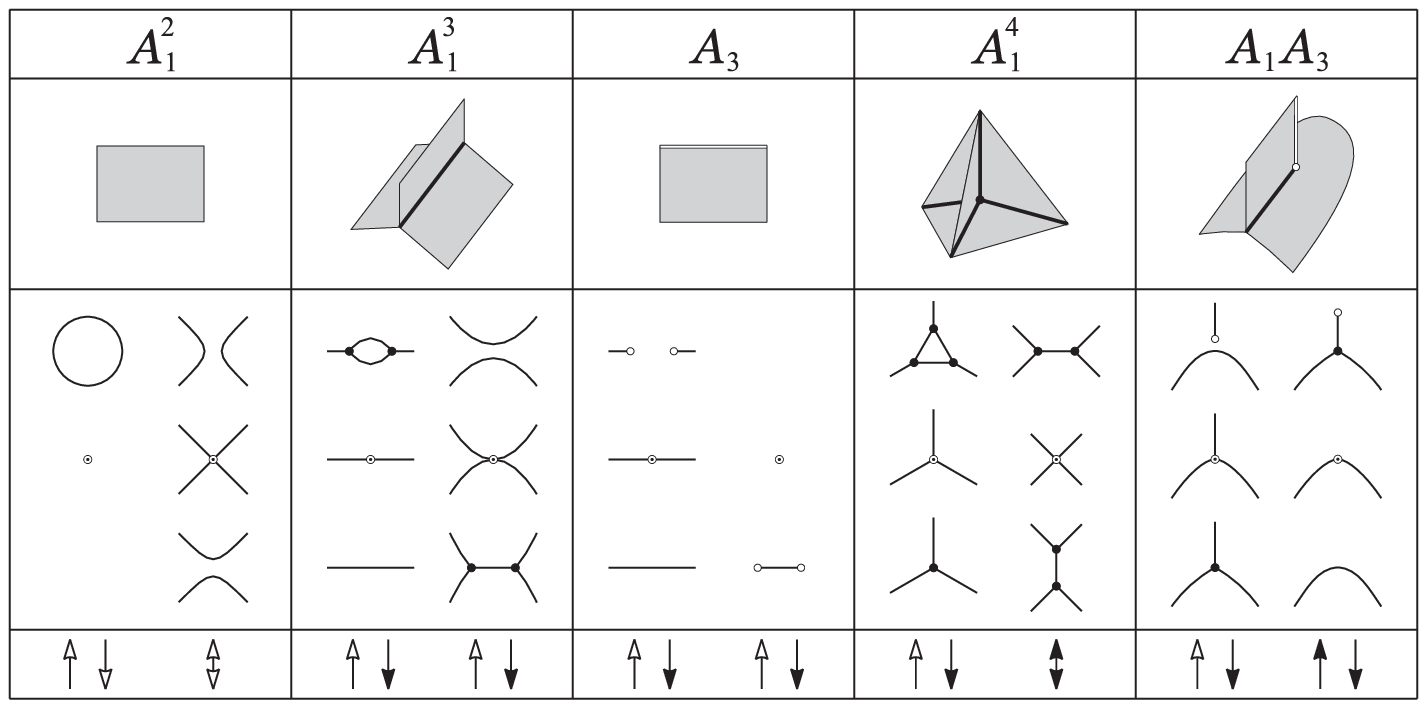}
\caption{Singularities of world shocks and perestroikas of instant
shocks in plane} \label{f2}
\end{center}
\end{figure*}

\begin{figure*}[t]
\begin{center}
\includegraphics{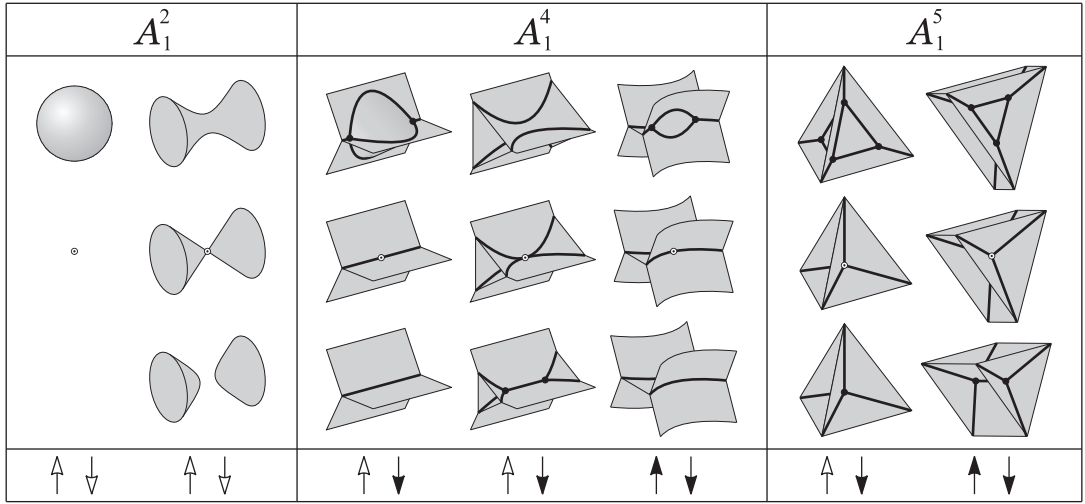}
\caption{Perestroikas of instant shocks in space (beginning)}
\label{f3a}
\end{center}
\end{figure*}

\begin{figure*}[t]
\begin{center}
\includegraphics{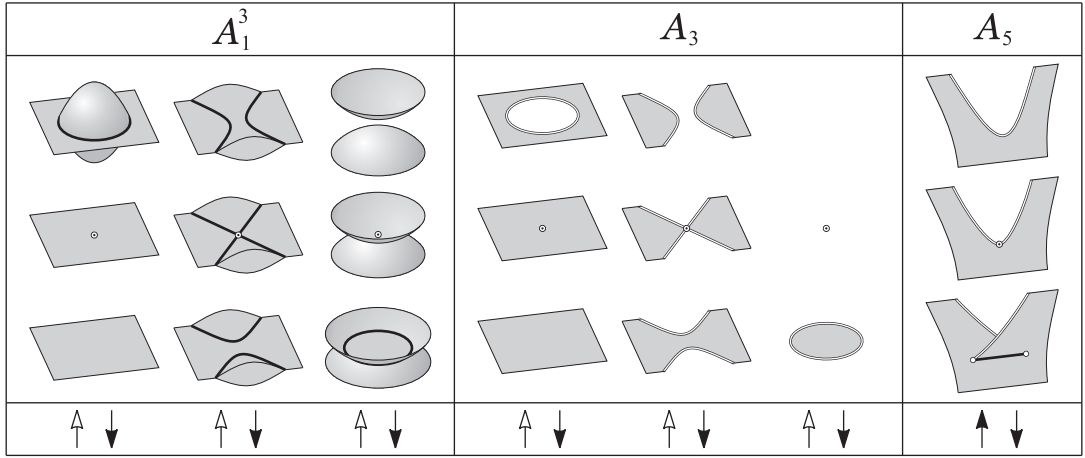}
\caption{Perestroikas of instant shocks in space (continuation)}
\label{f3b}
\end{center}
\end{figure*}

\begin{figure*}[t]
\begin{center}
\includegraphics{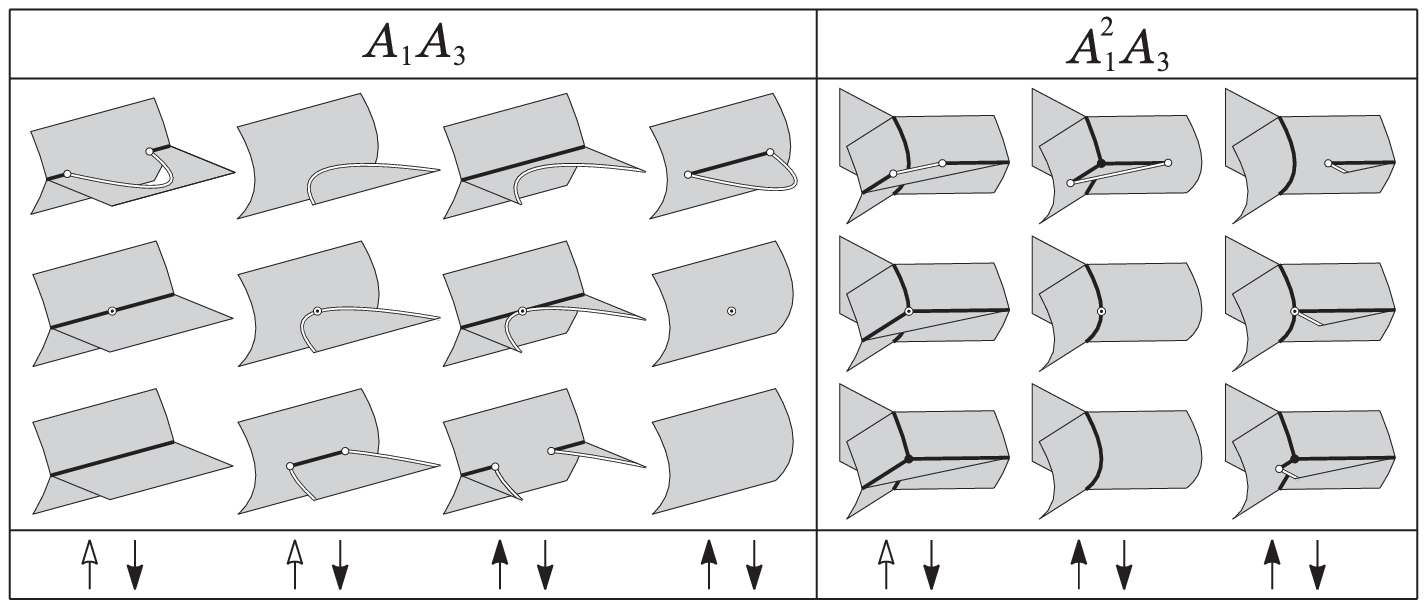}
\caption{Perestroikas of instant shocks in space (end)}
\label{f3c}
\end{center}
\end{figure*}

An example of a world shock in the well known case $d=1$ is shown
in Figure\,\ref{f1a}. The instant shocks are sets of isolated
points, the world shock is a curve in plane which can have triple
points and end-points. Any triple point sticks a pair of points of
an instant shock, any end-point generates its new point. These $2$
perestroikas exhaust all generic ones -- any other perestroika can
be killed by an arbitrarily small perturbation of the initial
condition.

Figure\,\ref{f2} shows the case $d=2$. A world shock is a surface
with singularities, all of them are shown in the second row. An
instant shock is a curve which can have triple points and
end-points -- the same singularities occur on world shocks of the
previous dimension $d=1$. At typical times there are no other
singularities, but at separate times they do occur. Namely, an
instant shock can experience perestroikas shown in the figure by
black arrows. (Up-down arrows show perestroikas which do not
change after the time inversion.) These $9$ perestroikas exhaust
all generic ones -- any other perestroika can be killed by an
arbitrarily small perturbation of the initial condition.

Figures \ref{f3a}, \ref{f3b}, and \ref{f3c} show the case $d=3$.
In this case it is impossible to draw singularities of world
shocks because they lie in $4$-dimensional space-time. An instant
shock is a surface which can have singularities which occur on
world shocks of the previous dimension $d=2$. Again, at typical
times there are no other singularities but at separate times they
occur. Namely, an instant shock can  experience perestroikas shown
in the figures by black arrows. These $26$ perestroikas exhaust
all generic ones -- any other perestroika can be killed by an
arbitrarily small perturbation of the initial condition.

In more detail, we consider a material $d$-dimensional medium
whose velocity is potential and described by the Burgers equation
with the potential force term:
\begin{equation}
\label{Burgers} \left\{
\begin{array}{l}
\vel_t + (\vel \cdot \nbl{}) \, \vel = -\nbl{U}
+ \visc \lpl{\vel}\\
\vel = \nbl{\potv}\\
\potv(\xb,0) = \pot_0(\xb)
\end{array}
\right.
\end{equation}
where $\xb \in \R^d$ is a point of the medium, $\vel(\xb,t)$ is
the velocity at the point $\xb$ at the time $t$, $\visc > 0$ is
the viscosity of the medium, $\nbl = (\partial_{x_1}, \dots,
\partial_{x_d})$ is the usual $\nabla$-operator in $\R^d$, and
$\lpl=\nabla \cdot \nabla$ is the Laplacian. The force potential
$U(\xb,t)$ and the initial condition $\pot_0(\xb)$ are assumed to
be {\it smooth}. (Everywhere in the present paper it means
infinite differentiability.) Let $\pot$ be the {\it limit
solution} as the viscosity vanishes. It has the following minimum
representation:
\begin{equation}
\label{F} \pot(\xb,t)=\lim_{\visc \to 0} \potv(\xb,t)= \min_\ub
F(\ub,\xb,t)
\end{equation}
where $F$ is a smooth function which, in the mechanics
terminology, is just the action defined by the Lagrangian
$$
\Lag(\dot{\xb},\xb,t)=\frac12 |\dot{\xb}|^2 - U(\xb,t)
$$
along its extremal trajectory ending at the point $(\xb,t)$ with
the velocity $\ub$. This fact is well known in the unforced case
$U \equiv 0$. If $U \not \equiv 0$ then it can be proved using the
so-called viscosity solutions of the Hamilton--Jacobi equation
described, for example, in \cite{CEL84}. In the unforced case
$\Lag=|\dot{\xb}|^2/2$
$$
F(\ub,\xb,t)=\pot_0(\xb-\ub t)+\frac12 |\ub|^2 t
$$
and we can rewrite our minimum representation as
\begin{equation}
\label{W}
\begin{array}{c}
{\displaystyle \pot(\xb,t)= \min_\ab \act(\ab,\xb,t)}\\[5pt]
{\displaystyle \act(\ab,\xb,t)= \pot_0(\ab) +
\frac{|\ab-\xb|^2}{2t}}
\end{array}
\end{equation}
where $\ab=\xb-\ub t$ and $\xb$ are the Lagrangian and Eulerian
coordinates of a particle moving with the velocity $\ub$. This
formula can be got explicitly without the viscosity solution
theory. We shall describe the both approaches in
Section\,\ref{vsHJ}.

The action $F$ is smooth but, nevertheless, the limit solution
$\pot$ is continuous only and can be non-differentiable at some
points. So, the velocity $\vel$ can have discontinuities. For
example, it can happen at a point $(\xb_\ast,t_\ast)$ such that
the function $F(\ub,\xb_\ast,t_\ast)$ of $\ub$ attains its minimum
for two different values $\ub_1 \ne \ub_2$:
$$
F(\ub_1,\xb_\ast,t_\ast) = F(\ub_2,\xb_\ast,t_\ast) \le
F(\ub,\xb_\ast,t_\ast)
$$
for any $\ub$.

\begin{figure}[h]
\begin{center}
\includegraphics{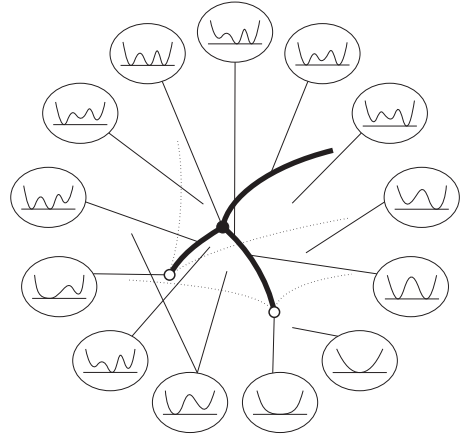}
\caption{Competition of minima generates shocks and caustics}
\label{bd}
\end{center}
\end{figure}

The points where the limit solution is not smooth form the {\it
world shock} lying in space-time and the {\it instant shocks}
which are in space and vary in time. In other words, the instant
shocks are the sections of the world shock with isochrones
$t=\cnst$. Figure\,\ref{bd} shows how competition of minima of the
functions $F(\ub,\xb,t)$ generates a world shock in
one-dimensional case $d=1$. Namely, at the points
$(\xb_\ast,t_\ast)$ of the world shock the corresponding function
$F(\ub,\xb_\ast,t_\ast)$ of $\ub$ has at least two global minima
or its minimum is degenerate. So, the world shock is the {\it
Maxwell stratum} of the family of functions of $\ub$ depending on
$\xb$ and $t$. The caustics are shown by dots -- at any their
point the corresponding function has a degenerate critical value
which is not minimal for a typical point of the caustic.

We describe local singularities and perestroikas of instant shocks
in physically interesting cases of dimension $d \le 3$ if the
initial condition is generic. Formally, it means that our
description works only for some open dense subset in the space of
all smooth initial conditions. It is assumed that this space is
provided with the so-called {\it Whitney topology} -- its
definition can be found, for example, in \cite{AGVI} (2.3,
Remark\,2). In other words, we are interested only in
singularities and perestroikas which are stable with respect to
any sufficiently small perturbations of the initial condition.
Other singularities and perestroikas can be killed by its
arbitrarily small perturbation.

It turns out that such a generic instant shock is a hypersurface
which at typical times has singularities from a finite (up to
smooth diffeomorphisms) list.  At these times the instant shock
does not locally change with respect to smooth diffeomorphisms.
But at separate times it can experience local perestroikas -- the
instant shocks at the perestroika time and at close times are not
diffeomorphic to each other. The instant shock at the perestroika
time itself has a singularity which does not belong to the finite
list.

We get a finite topological classification of generic perestroikas
of instant shocks.  Their smooth classification proves to be
infinite -- continuous invariants or {\it moduli} appear. The
finite topological classification described in detail in
Section\,\ref{res} is the first main result of the present paper
(Theorem\,\ref{mr1}). The second main result is the shock
determinator described in Section\,\ref{det} (Theorems \ref{sing},
\ref{pern}, and \ref{pers}) which allows to define what
singularity or perestroika is occurring at a given point of
space-time for a concrete initial condition of the unforced
Burgers equation.

The general plan of the classification is the following.


\subsection{Step 1: Minimum representations for limit solutions}

Firstly, we describe in Section\,\ref{vsHJ} how to get the minimum
representations of limit solutions. The representation (\ref{W})
can be obtained finding explicitly the potential solution of the
unforced Burgers equation with the help of the so-called
Hopf--Cole transformation which was published by Forsyth back in
1906 \cite{For06}. This is a well known way, but we describe it.

The representation (\ref{F}) in the forced case is an automatic
corollary of the theory of generalized or viscosity solutions of a
Hamilton--Jacobi equation constructed in \cite{Kr66}, \cite{Kr67},
\cite{CEL84}, we brief this approach in Section\,\ref{vsHJ} as
well.

As a matter of fact, our description of perestroikas works for the
shocks of viscosity or generalized solutions of the
Hamilton--Jacobi equation with a strictly convex (up to momenta)
smooth Hamiltonian \cite{Bog99}. In this case there exists a
minimum representation which is analogous to (\ref{F}) and
published, for example, in \cite{ES84}. This representation shows
that the viscosity solution in the convex case is the solution of
the classical variational problem with a free left end:
$$
\pot(\xb,t) =\min_{\gamma(t)=\xb} \left\{ \pot_0(\gamma(0)) +
\int\limits_0^{t} \Lag(\dot{\gamma} (\tau), \gamma (\tau), \tau)
\,\df\tau \right\}
$$
where $\gamma: [0,t] \to \R^d$ is a differentiable path in space
and the Lagrangian is the Legendre transformation of the
Hamiltonian:
$$
\Lag(\dot{\xb},\xb,t) = \max_\pb{\left\{\pb \dot{\xb} -
\Ham(\pb,\xb,t) \right\}}.
$$
The Hamiltonian and the Lagrangian are both strictly convex down
with respect to momenta. In other words, the matrices of their
second derivatives must be positive definite:
$$
\left\|\Ham_{\pb\pb} (\pb,\xb,t)\right\| > 0, \quad
\left\|\Lag_{\dot{\xb}\dot{\xb}} (\dot{\xb},\xb,t)\right\| > 0.
$$

Therefore, our classification describes singularities and
perestroikas of solutions of the above variational problem. (The
Lagrangian $\Lag$ is assumed to be fixed but the initial condition
$\pot_0$ is generic.) In fact, everything is vice versa: at first
we classify singularities and perestroikas in the variational
problem with a convex Lagrangian, and after that we apply this
classification to the limit solutions taking into account the
minimum representations.

\subsection{Step 2: Perestroikas of shocks for generic $F$}

Secondly, we forget that $F$ is the action and consider it just as
a generic smooth function. In Section\,\ref{minfuns} we describe,
using methods of Singularity Theory, all singularities and
perestroikas of shocks for $d \le 3$ which are stable with respect
to small perturbation of the function $F$ itself. Other
singularities and perestroikas can be killed by an arbitrarily
small perturbation of the function $F$.

Singularities of instant shocks at typical times and singularities
of world shocks were described in \cite{Br77}. Perestroikas of
instant shocks were published for the first time in \cite{Bog89}.
The both classification lists are proved to be finite for
appropriate equivalence relations: smooth for singularities and
topological for perestroikas.


But can the action $F$ be considered as generic? Indeed, $F$ is
not an arbitrary function, it is defined by the initial condition
and the force potential. For example, in the unforced case there
is the formula (\ref{W}) and the action does not allow any
perturbation as a function of Lagrangian coordinates $\ab$,
Eulerian coordinates $\xb$, and time $t$. Only the initial
condition $\pot_0$ is a generic smooth function of the Lagrangian
coordinates $\ab$ -- the other part of the formula (\ref{W}) is a
fixed function. Does this fact give us new singularities and
perestroikas of shocks which are stable with respect to small
perturbations of the initial condition only but not the action
itself? Maybe, there exist singularities or perestroikas of shocks
which cannot be killed by an arbitrarily small perturbation of the
initial condition but vanish after such a perturbation of the
action?

It turns out that the answers to these questions are negative -- a
general theorem implying this fact is published in \cite{Bog95}.
Therefore, all singularities and perestroikas of instant shocks
with generic initial conditions belong to our lists.

\subsection{Step 3: Topological restrictions forbidding some perestroikas}

On the other hand, do there exist singularities or perestroikas
from the list of the step\,2 which cannot be realized by shocks
when $F$ is really the action? It turns out that such
singularities do not exist but there are a few perestroikas
forbidden for any force potential $U$. Apparently, for the first
time this fact has been observed by S.\,N.\,Gurbatov and
A.\,I.\,Saichev. Firstly, an isochrone $t=\cnst$ cannot be tangent
to a world shock at its smooth point -- this restriction forbids
all perestroikas $A_1^2$. Secondly, a few perestroikas are
irreversible -- they can occur in only one direction of time: for
example, a triangle from shocks can vanish during the perestroika
$A_1^3$ but cannot appear. A universal topological restriction is
the following:

{\it A local shock after a perestroika is contractible
(homotopically equivalent to a point).} For example, the triangle
from shocks is homotopically a circle but not a point, hence it
cannot appear.

This topological restriction was explicitly checked in
\cite{Bog89} only for the perestroikas from the list of the
step\,2. Namely, the forbidden perestroikas cannot be realized by
any action in view of algebraic conditions described in
Section\,\ref{shkrec} for the unforced case. In
Section\,\ref{toprstr} we prove this restriction for an
essentially wider class of perestroikas in all dimensions. As a
matter of fact, it is almost equivalent to another topological
restriction suggested and proved in \cite{Bar90} for all
dimensions and not only for generic perestroikas. Who knows why
nobody noticed this simple equivalence for a dozen of years?

\subsection{Step 4: Realization of not forbidden perestroikas}

At last we show that all perestroikas not forbidden by the
topological restriction can be realized by instant shocks of the
Burgers equation. Moreover, Theorems of Section\,\ref{det} allow,
knowing the initial condition, to define what singularity or
perestroika from our list is realized by the instant shock at a
given point of space-time. In particular, they give all genericity
conditions for the initial condition which guarantee that the
instant shock has only described singularities and experiences
only described perestroikas. These theorems are proved in
Section\,\ref{shkrec}.

\section{The Classification}

\label{res}

\subsection{Simple minima}

If in proper local smooth (curvilinear, in general!) coordinates
a function on $\R^d$ has the form
$$\ui_1^{2\nk} + \ui_2^2 + \dots + \ui_d^2 + \cnst$$
then it attains a local minimum at the point $\ub=0$. In
Singularity Theory such minimum has the name $A_{2\nk-1}$. The
index $\mu=2\nk-1$ is called the {\it multiplicity} of the
minimum. The number $\kappa=\mu-1=2\nk-2$ is called the {\it
codimension} of the minimum. Any non-degenerate minimum of a
smooth function is $A_1$ -- this is the well known Morse lemma.
The minima $A_3$, $A_5$,... are degenerate, the ranks of their
second differentials are equal to $d-1$. So, it is easy to find a
minimum which does not belong to the series $A_{2\nk-1}$ -- the
simplest one is:
$$
X_9 : \; u^4 + 2\,C\,u^2 v^2 + v^4 + \cnst, \quad C > -1, \; C \ne 1,
$$
where $C$ is a modulus or a continuous invariant up to smooth
changes of the variables. (If $C=1$ then $0$ is a minimum too but
its name is not $X_9$.)

\begin{figure}[h]
\begin{center}
\includegraphics{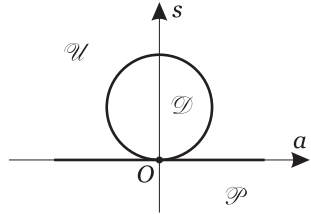}
\caption{Kernel circle and line divide kernel plane of minimum
$A_3$ into three domains} \label{k}
\end{center}
\end{figure}

Any smooth function with a minimum $A_3$ in $\R^d$ can be
represented in some orthonormal coordinates $(a,\bb)$ with origin
at the minimum point as
\begin{equation}
\label{ra3}
\begin{array}{c}
K a^4 + \omega(\vb) + O_{>4}(a,\vb) + \cnst,\\[5pt]
\omega(\vb) = \omega_1 \vi_1^2 + \dots + \omega_{d-1}
\vi_{d-1}^2,\\[5pt]
\vb = \bb - \Cb a^2
\end{array}
\end{equation}
where $\bb,\vb,\Cb \in \R^{d-1}$, $K, \omega_1, \dots,
\omega_{d-1} > 0$, and $O_{>4}$ denotes terms of order greater
than 4 if $\deg a=1$ and $\deg \vb = 2$. The straight line $\bb=0$
is invariant and called {\it kernel}. If $\Cb \ne 0$ then the
points
$$
(a, s \, \Cb / |\Cb|), \quad (a,s) \in \R^2
$$
form the invariant {\it kernel} plane, and the equation
$$|\Cb|(a^2+s^2)-s=0$$
defines the {\it kernel} circle which lies in the kernel plane and
is tangent to the kernel line $s=0$. The kernel circle and line
divide the kernel plane into three domains: the disk
$\mathcal{D}$, the domain $\mathcal{U}$, and the semiplane
$\mathcal{P}$ shown in Figure\,\ref{k}.

Any smooth function with a minimum $A_5$ in $\R^d$ can be
represented in some orthonormal coordinates $(a,\bb)$ with origin
at the minimum point as
\begin{equation}
\label{ra5}
\begin{array}{c}
K a^6 + \omega(\vb) + O_{>6}(a) + O_{>8}(a,\vb) + \cnst,\\[5pt]
\omega(\vb) = \omega_1 \vi_1^2 + \dots + \omega_{d-1}
\vi_{d-1}^2,\\[5pt]
\vb = \bb - \Cb a^2 - \mathbf{D}a^3 - \mathbf{E}a^4
\end{array}
\end{equation}
where $\bb,\vb,\Cb,\mathbf{D},\mathbf{E} \in \R^{d-1}$, $K,
\omega_1, \dots,\omega_{d-1} > 0$, and $O_{>n}$ denotes terms of
order greater than $n$ if $\deg a=1$ and $\deg \vb = 4$.

\subsection{Classification of shock points}

If a point $(\xb_\ast,t_\ast)$ lies on a shock then the function
$F(\cdot,\xb_\ast,t_\ast)$ attains its minimal value at a few
points or at a degenerate one. The names of these minima can be
written as one word like $A_1^2A_3$ which means: ``two
(non-degenerate) minima $A_1$ and one (degenerate) minimum
$A_3$''. Such word is the name of the point $(\xb_\ast,t_\ast)$.

For example, a typical point of a generic shock is $A_1^2$. In its
neighborhood the shock is a smooth hypersurface which acquires
singularities at the points $A_1^3$, $A_3$,... The points $A_1$
are out of the shock at all.

\subsection{Shocks on line: $d=1$}

\begin{figure}[h]
\begin{center}
\includegraphics{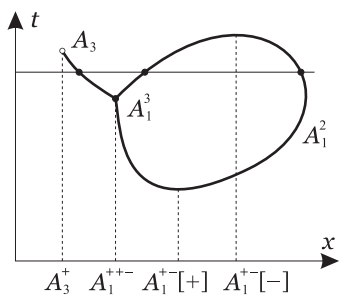}
\caption{Impossible perestroikas of instant shocks on line}
\label{f1b}
\end{center}
\end{figure}

If the action $F$ is just an arbitrary generic smooth function,
the world shock is a curve in plane which, besides smooth points
$A_1^2$, can have triple points $A_1^3$ and end-points $A_3$.  A
generic curve with such prescribed singularities could generate
$6$ different perestroikas of its sections which are just sets of
isolated points. All these possibilities -- $A_1^2[\pm]$ (a couple
of points appears or vanishes), $A_1^{+\pm-}$, (a point doubles or
two points stick), $A_3^\pm$ (a point vanishes or appears) -- are
shown in Figures \ref{f1a} and \ref{f1b}.

Which of these perestroikas are realized by ``real'' shocks,
i.\,e. the shocks of limit potential solutions of the Burgers
equation? According to the well known entropy condition, only two
$A_1^{+--}$ (two points stick) and $A_3^-$ (a point appears) can
happen -- they are shown in Figure\,\ref{f1a}. The $4$ forbidden
perestroikas are shown in Figure\,\ref{f1b}. They can be explained
by our topological restriction as well: the empty set and a couple
of points are not contractible.

\subsection{Shocks in plane: $d=2$}

Let the action $F$ be an arbitrary generic smooth function. Then
the instant shock at a typical time is a curve which, besides
smooth points $A_1^2$, can have triple points $A_1^3$ and
end-points $A_3$ -- the same singularities occur on world shocks
of the previous dimension $d=1$. The world shock is a surface
which, besides smooth points $A_1^2$, can have the following
singularities: $A_1^3$, $A_3$, $A_1^4$, and $A_1A_3$ shown in
Figure\,\ref{f2}. The singularities $A_1^4$, and $A_1A_3$ are new
in comparison with the previous dimension and occur at isolated
points of the world shock. The singularities $A_1^3$ and $A_3$
form curves -- triple lines and edges.

A generic surface with such prescribed singularities generates
$18$ different perestroikas of its sections with isochrones -- all
of them are shown in Figure\,\ref{f2} by white and black arrows.
The perestroika $A_1^2$ occurs when an isochrone is tangent to the
surface at a smooth point, the perestroikas $A_1^3$ and $A_3$
happen if an isochrone is tangent to a triple line or to an edge.
But only $9$ perestroikas are not forbidden  by the topological
restriction and can be realized by ``real'' shocks. They are shown
by black arrows. (White ones show the forbidden perestroikas.)

\subsection{Shocks in space: $d=3$}

Again, let the action $F$ be an arbitrary generic smooth function.
Then the instant shock at a typical time is a surface which,
besides smooth points $A_1^2$, can have triple lines $A_1^3$,
edges $A_3$, and separate singularities $A_1^4$ and $A_1A_3$ --
the same singularities occur on world shocks of the previous
dimension $d=2$. A world shock, besides smooth points $A_1^2$, can
have the following singularities: $A_1^3$, $A_3$, $A_1^4$,
$A_1A_3$, $A_1^5$, $A_1^2A_3$, and $A_5$. The singularities
$A_1^5$, $A_1^2A_3$, and $A_5$ are new in comparison with the
previous dimension and occur at isolated points of the world
shock. The singularities $A_1^4$ and $A_1A_3$ form curves and the
singularities $A_1^3$ and $A_3$ form surfaces.

In this case there are $42$ different perestroikas -- all of them
are shown in Figures \ref{f3a}, \ref{f3b}, and \ref{f3c} by white
and black arrows. But only $26$ of them are not forbidden  by the
topological restriction and can be realized by  ``real'' shocks.
They are shown by black arrows. (White ones show the forbidden
perestroikas.)

\subsection{Normal forms}

Formulas describing singularities of shocks are given in
Table\,\ref{t}. The functions from the second column of
Table\,\ref{t} are called {\it normal forms of a limit solution}
in local smooth coordinates $(\alpha, \beta, \dots)$ with respect
to adding a smooth function. The set of points, where a normal
form is not smooth, describes a world shock in a neighborhood of a
point with the indicated name. Its coordinates are $\alpha = \beta
= \dots =0$. In general, the coordinates $(\alpha, \beta, \dots)$
are curvelinear and mix space coordinates $\xb$ and time $t$.
Besides, these normal forms describe limit solutions and instant
shocks at typical times. In this case $(\alpha, \beta, \dots)$ are
local smooth coordinates in space.

\begin{table*}[t]
\begin{center}
\begin{tabular}{|c|c|c|c|c|}
\hline
&&\multicolumn{3}{c|}{}\\[-9pt]
Name & Limit Solution & \multicolumn{3}{c|}{
\begin{tabular}{c}
Time\\
$\phantom{d=1}$
\end{tabular}
}
\\[-13pt]
\cline{3-5}
&&&&\\[-9pt]
&& $d=1$ & $d=2$ & $d=3$
\\
\hline
&&&&\\[-10pt]
$A_1^2$ & $\min \{\alpha,0\}$ & $\pm(\alpha + \beta^2)$ &
\begin{tabular}{c}
$\pm(\alpha + \beta^2 +\gamma^2)$\\
$\alpha + \beta^2 -\gamma^2$
\end{tabular}
& \begin{tabular}{c}
$\pm(\alpha + \beta^2 +\gamma^2 +\delta^2)$\\
$\pm(\alpha + \beta^2 +\gamma^2 -\delta^2)$
\end{tabular}
\\
\hline
&&&&\\[-10pt]
$A_1^3$ & $\min \{\alpha,\beta,0\}$ & $\pm(\alpha + \beta)$ &
\begin{tabular}{c}
$\pm(\alpha + \beta +\gamma^2)$\\
$\pm(\alpha + \beta -\gamma^2)$
\end{tabular}
&
\begin{tabular}{c}
$\pm(\alpha + \beta +\gamma^2 +\delta^2)$\\
$\pm(\alpha + \beta +\gamma^2 -\delta^2)$\\
$\pm(\alpha + \beta -\gamma^2 -\delta^2)$
\end{tabular}
\\
\hline
&&&&\\[-10pt]
$A_1^4$ & $\min \{\alpha,\beta,\gamma,0\}$ &&
\begin{tabular}{c}
$\pm(\alpha + \beta +\gamma)$\\
$\alpha + \beta -\gamma$
\end{tabular}
&
\begin{tabular}{c}
$\pm(\alpha + \beta +\gamma +\delta^2)$\\
$\pm(\alpha + \beta +\gamma -\delta^2)$\\
$\pm(\alpha + \beta -\gamma +\delta^2)$
\end{tabular}
\\
\hline
&&&&\\[-10pt]
$A_1^5$ & $\min \{\alpha,\beta,\gamma,\delta,0\}$ &&&
\begin{tabular}{c}
$\pm(\alpha + \beta +\gamma +\delta)$\\
$\pm(\alpha + \beta +\gamma -\delta)$
\end{tabular}
\\
\hline
&&&&\\[-10pt]
$A_3$ & $\min\limits_\ui (\ui^4 + \alpha \ui^2 + \beta \ui)$ &
$\pm\alpha$ &
\begin{tabular}{c}
$\pm(\alpha+\gamma^2)$\\
$\pm(\alpha-\gamma^2)$
\end{tabular}
&
\begin{tabular}{c}
$\pm(\alpha +\gamma^2 +\delta^2)$\\
$\pm(\alpha +\gamma^2 -\delta^2)$\\
$\pm(\alpha -\gamma^2 -\delta^2)$
\end{tabular}
\\
\hline
&&&&\\[-10pt]
$A_1A_3$ & $\min \left\{\alpha, \min\limits_\ui (\ui^4 + \beta
\ui^2 + \gamma \ui) \right\}$ &&
\begin{tabular}{c}
$\pm(\alpha +\beta)$\\
$\pm(\alpha -\beta)$
\end{tabular}
&
\begin{tabular}{c}
$\pm(\alpha +\beta +\delta^2)$\\
$\pm(\alpha +\beta -\delta^2)$\\
$\pm(\alpha -\beta -\delta^2)$\\
$\pm(\alpha -\beta +\delta^2)$
\end{tabular}
\\
\hline
&&&&\\[-10pt]
$A_1^2A_3$ & $\min \left\{\alpha,\beta, \min\limits_\ui (\ui^4 +
\gamma \ui^2 + \delta \ui) \right\}$ &&&
\begin{tabular}{c}
$\pm(\alpha +\beta +\gamma)$\\
$\pm(\alpha +\beta -\gamma)$\\
$\pm(\alpha -\beta +\gamma)$
\end{tabular}
\\
\hline
&&&&\\[-10pt]
$A_5$ & $\min\limits_\ui (\ui^6 + \alpha \ui^4 + \beta \ui^3 +
\gamma \ui^2 + \delta \ui)$ &&&
$\pm\alpha$\\
\hline
\end{tabular}
\caption{Normal forms of singularities and perestroikas of shocks}
\label{t}
\end{center}
\end{table*}

Perestroikas of instant shocks are described by the {\it normal
forms of time} from Table\,\ref{t}. The world shock is described
by the second column, the time coordinate is
$$t(\alpha,\beta,\dots) = t_\ast + \tau(\alpha, \beta, \dots)$$
where $t_\ast$ is the perestroika's instant and $\tau$ is a normal
form from the following three columns which correspond to the
dimensions $d=1,2,3$ respectively. An empty cell of Table\,\ref{t}
means that the corresponding singularity of a world shock does not
occur in the shown dimension.

The perestroikas for $d=1$ from Table\,\ref{t} are shown in
Figures \ref{f1a} and \ref{f1b}. Plus or minus denoted by the sign
``$\pm$'' in the normal form of time from the table corresponds to
the same sign in perestroikas' names $A_1^{+-}[\pm]$,
$A_1^{+\pm-}$, and $A_3^\pm$ in the figures.

The perestroikas for $d=2,3$ from Table\,\ref{t} are shown in
Figures \ref{f2}, \ref{f3a}, \ref{f3b}, and \ref{f3c} by arrows.
Namely, each perestroika is shown as the three sections of the
world shock with isochrones
$$
\begin{array}{l}
t(\alpha, \beta, \dots) = t_- < t_\ast,\\
t(\alpha, \beta, \dots) = t_\ast,\\
t(\alpha, \beta, \dots) = t_+ > t_\ast.
\end{array}
$$
An arrow shows the time direction from $t_-$ to $t_+$. The
$t_\pm$-sections have only generic singularities, but one point of
the $t_\ast$-section is not generic -- it is the point of the
perestroika, its coordinates are $\alpha=\beta=\dots=0$. The
perestroikas which do not change after the time inversion are
shown by up-down arrows $\updownarrow$ in the figures, their
functions $\tau$ from Table\,\ref{t} do not contain $\pm$.  The
sign ``$\pm$'' in Table\,\ref{t} means that the perestroika does
change after the time inversion, $\uparrow$ in the figures means
$+$, $\downarrow$ means $-$.

Unfortunately, if we considered only smooth coordinates it would
be impossible to get such simple normal forms of time -- there
would be quite a lot of continuous invariants. In other words,
smooth classification of perestroikas of instant shocks is
infinite. So, in order to describe perestroikas of shocks we allow
the coordinates $(\alpha, \beta, \dots)$ to be sectionally
continuously differentiable but, nevertheless, they must be smooth
outside the world shock. Therefore, the first derivatives
$\partial (\alpha, \beta, \dots)/\partial (\xb,t)$ can have jumps
on the world shock itself.

The first main result of the present paper is

\medskip
\begin{theorem}
\label{mr1} Let $d \le 3$, $T>0$, and the initial condition be a
generic smooth function. Then any perestroika of the instant shock
happening before the time $T$ is described in some coordinates
$(\alpha,\beta,\dots)$ by a pair of normal forms from
Table\,\ref{t} -- the limit solution and time. The coordinates
$(\alpha, \beta, \dots)$ are sectionally continuously
differentiable and smooth outside the world shock of the limit
solution.

But not all perestroikas from Table\,\ref{t} occur -- for real
ones {\bf the local instant shock after the perestroika must be
contractible}. Such perestroikas are indicated in Figure\,\ref{f2}
($d=2$) and Figures \ref{f3a}, \ref{f3b}, \ref{f3c} ($d=3$) by
black arrows, the well known case $d=1$ is shown in
Figure\,\ref{f1a}.
\end{theorem}

\medskip

\begin{remark}
The term ``generic'' means ``from some open dense subset''. It is
assumed that the space of all smooth initial conditions is
provided with the so-called {\it Whitney topology} -- its
definition can be found, for example, in \cite{AGVI} (2.3,
Remark\,2). If $T=\infty$ then the theorem is true only for some
dense subset of the initial conditions which maybe is not open.
\end{remark}

\begin{remark}
Another way to interpret our classification of perestroikas is to
use smooth coordinates but another equivalence for limit
solutions. Namely, instead of adding a smooth function we allow a
smooth diffeomorphism of the graphs of the limit solution and its
normal form at the same time. All details are discussed in
\ref{pmf}.
\end{remark}

\subsection{Smoothness of $F$}

In the present paper we require infinite smoothness of the action
$F$ which can be guaranteed by infinite smoothness of the initial
condition $\pot_0$ and the force potential $U$. Probably,
Theorem\,\ref{mr1} is correct when we can guarantee $F \in
C^{d+3}$ with respect to all variables ($\ub$, $\xb$, and $t$).
Then the coordinates $(\alpha, \beta, \dots)$ from
Theorem\,\ref{mr1} will be sectionally continuously differentiable
as before but only continuously differentiable outside the world
shock of the limit solution. The crucial point in this question is
Lemma\,\ref{fams} -- the change $h$ must be continuously
differentiable at least.

\section{The Shock Determinator}

\label{det}

The worth of Theorem\,\ref{mr1} for applications is problematic.
How can one find out whether or not a given initial condition is
generic? What does really happen at some point of space-time if the
initial condition is fixed and cannot be perturbed? Briefly
speaking, how is it possible to apply Theorem\,\ref{mr1} in a
concrete situation?

The answer for the unforced Burgers equation is given by Theorems
\ref{sing}, \ref{pern} and \ref{pers}. Theorem\,\ref{sing}
recognizes the singularity of the instant shock at a given point
of space-time. If none of its conditions is satisfied then we can
try to apply Theorems \ref{pern} and \ref{pers} which recognize
the perestroika of the instant shock. If they do not work too, it
means that something non-generic occurs at our point. Namely, the
singularity or perestroika is not stable with respect to small
perturbations of the initial condition -- it is not generic and,
therefore, is not described by Theorem\,\ref{mr1}.

For the unforced Burgers equation there is the representation
(\ref{W}). Let $(\xb_\ast,t_\ast)$ be a fixed point of space-time
and
$$
\sact(\cdot) = \act(\cdot,\xb_\ast,t_\ast).
$$
The information about global minima of the function $\sact$ is
written like $\sact \in A_1^2A_3$ -- this example means that
$\sact$ has three global minima, two of them are $A_1$ and the
other one is $A_3$.

\medskip
\begin{theorem}
\label{sing} {\sc (Singularities)} If $d=1,2,3$ and $\sact \in
A_1$, then $\xb_\ast$ lies outside the instant shock at the time
$t_\ast$.

If $d=1,2,3$ and $\sact \in A_1^2$, then $\xb_\ast$ is a smooth
point $A_1^2$ of the instant shock at the time $t_\ast$.

If $d=2,3$, $\sact \in A_1^3$, and the minimum points do not
belong to the same straight line, then $\xb_\ast$ is a singularity
$A_1^3$ of the instant shock at the time $t_\ast$.

If $d=3$, $\sact \in A_1^4$, and the minimum points do not belong
to the same plane, then $\xb_\ast$ is a singularity $A_1^4$ of the
instant shock at the time $t_\ast$.

If $d=2,3$, $\sact \in A_3$ such that $\Cb \ne 0$, then $\xb_\ast$
is a singularity $A_3$ of the instant shock at the time $t_\ast$.

If $d=3$, $\sact \in A_1 A_3$, at the minimum $A_3$ $\Cb \ne 0$,
and the kernel plane does not contain the minimum point $A_1$,
then $\xb_\ast$ is a singularity $A_1A_3$ of the instant shock at
the time $t_\ast$.
\end{theorem}
\medskip

In order to recognize perestroikas we use their extended names and
signatures.

The {\it extended name} of a perestroika codes its name and the
linear part of the normal form of time from Table\,\ref{t}.  For
example, the normal form $-\alpha+\beta-\gamma$ in the case
$A_1^2A_3$ gives the extended name $A_1^{+-}A_3^-$. Namely, the
two first signs corresponding to the two minima $A_1$ (see the
normal form of limit solution) give $A_1^{+-}$ -- the order is
irrelevant, the third sign corresponding to the minimum $A_3$
gives $A_3^-$. The cases $A_1^\nm$ are more complicated: the
extended name shows not only the signs of the coefficients of the
linear part but the opposite sign of their sum as well. For
example, the normal form $-\alpha-\beta-\gamma+\delta$ in the case
$A_1^5$ gives the extended name $A_1^{++---}$ -- the signs
themselves give $---+$ and $+$ is the sign of the number
$-(-1-1-1+1)=2$, the order is irrelevant again. A general
definition will be given in \ref{pmf}.

The {\it signature} of a perestroika is just the signature of the
quadratic part of the normal form of time from Table\,\ref{t} --
the quantities of the positive and negative squares respectively.
The signature is written as $[p,q]$, the trivial ones $[0,0]$ are
omitted. For example, the normal form
$-\alpha-\beta^2-\gamma^2+\delta^2$ in the case $A_1^2$ gives the
signature $[1,2]$ -- the order is irrelevant as usual. The
corresponding perestroika is recognized as $A_1^{+-}[1,2]$.

The extended names and signatures of the perestroikas from
Figure\,\ref{f2} and Figures\,\ref{f3a},\ref{f3b},\ref{f3c} are
shown in Tables \ref{map2} and \ref{map3} respectively. In other
words, Tables \ref{map2} and \ref{map3} can be considered as
``maps'' of the figures of perestroikas. Namely, the top line of
each cell of the tables indicates the extended name and signature
of the perestroika $\uparrow$ from the corresponding cell of the
figures, the bottom line corresponds to the perestroika
$\downarrow$.

\begin{table*}[t]
\begin{center}
\begin{tabular}{|c|c|c|c|c|c|c|c|c|c|}
\hline \multicolumn{2}{|c|}{} & \multicolumn{2}{c}{} &
\multicolumn{2}{|c|}{} & \multicolumn{2}{c}{} &
\multicolumn{2}{|c|}{}
\\[-10pt]
\multicolumn{2}{|c|}{$A_1^2$} & \multicolumn{2}{c}{$A_1^3$} &
\multicolumn{2}{|c|}{$A_3$} & \multicolumn{2}{c}{$A_1^4$} &
\multicolumn{2}{|c|}{$A_1A_3$}
\\
\hline
&&&&&&&&&\\[-10pt]
$A_1^{+-}[++]$ & $A_1^{+-}[+-]$ & $A_1^{++-}[+]$ & $A_1^{++-}[-]$
& $A_3^+[+]$ & $A_3^+[-]$ & $A_1^{+++-}$ & $A_1^{++--}$ &
$A_1^+A_3^+$ & $A_1^+A_3^-$
\\
&&&&&&&&& \\[-10pt]
$A_1^{+-}[--]$ & $A_1^{+-}[+-]$ & $A_1^{+--}[-]$ & $A_1^{+--}[+]$
& $A_3^-[-]$ & $A_3^-[+]$ & $A_1^{+---}$ & $A_1^{++--}$ &
$A_1^-A_3^-$ & $A_1^-A_3^+$
\\
\hline
\end{tabular}
\caption{Map of perestroikas of shocks in plane} \label{map2}
\end{center}
\end{table*}

\begin{table*}[t]
\begin{center}
\begin{tabular}{|c|c|c|c|c|c|c|}
\hline \multicolumn{2}{|c|}{} & \multicolumn{3}{c}{} &
\multicolumn{2}{|c|}{}
\\[-10pt]
\multicolumn{2}{|c|}{$A_1^2$} & \multicolumn{3}{c}{$A_1^4$} &
\multicolumn{2}{|c|}{$A_1^5$}
\\
\hline
&&&&&&\\[-10pt]
$A_1^{+-}[+++]$ & $A_1^{+-}[++-]$ & $A_1^{+++-}[+]$ &
$A_1^{+++-}[-]$ & $A_1^{++--}[+]$ & $A_1^{++++-}$ & $A_1^{+++--}$
\\
&&&&&&\\[-10pt]
$A_1^{+-}[---]$ & $A_1^{+-}[+--]$ & $A_1^{+---}[-]$ &
$A_1^{+---}[+]$ & $A_1^{++--}[-]$ & $A_1^{+----}$   &
$A_1^{++---}$
\\
\hline
\multicolumn{7}{c}{}\\[-8pt]
\hline \multicolumn{3}{|c|}{} & \multicolumn{3}{c}{} &
\multicolumn{1}{|c|}{}
\\[-10pt]
\multicolumn{3}{|c|}{$A_1^3$} & \multicolumn{3}{c}{$A_3$} &
\multicolumn{1}{|c|}{$A_5$}
\\
\hline
&&&&&&\\[-10pt]
$A_1^{++-}[++]$ & $A_1^{++-}[+-]$ & $A_1^{++-}[--]$ & $A_3^+[++]$
& $A_3^+[+-]$ & $A_3^+[--]$ & $A_5^+$
\\
&&&&&&\\[-10pt]
$A_1^{+--}[--]$ & $A_1^{+--}[+-]$ & $A_1^{+--}[++]$ & $A_3^-[--]$
& $A_3^-[+-]$ & $A_3^-[++]$ & $A_5^-$
\\
\hline
\multicolumn{7}{c}{}\\[-8pt]
\hline
\multicolumn{4}{|c|}{} & \multicolumn{3}{c|}{}\\[-10pt]
\multicolumn{4}{|c|}{$A_1A_3$} &
\multicolumn{3}{c|}{$A_1^2A_3$}\\
\hline
&&&&&&\\[-10pt]
$A_1^+A_3^+[+]$ & $A_1^+A_3^+[-]$ & $A_1^+A_3^-[-]$ &
$A_1^+A_3^-[+]$ & $A_1^{++}A_3^+$ & $A_1^{++}A_3^-$ &
$A_1^{+-}A_3^+$
\\
&&&&&&\\[-10pt]
$A_1^-A_3^-[-]$ & $A_1^-A_3^-[+]$ & $A_1^-A_3^+[+]$ &
$A_1^-A_3^+[-]$ & $A_1^{--}A_3^-$ & $A_1^{--}A_3^+$ &
$A_1^{+-}A_3^-$
\\
\hline
\end{tabular}
\caption{Map of perestroikas of shocks in space} \label{map3}
\end{center}
\end{table*}

\medskip
\begin{theorem}
\label{pern} {\sc (Perestroikas -- Extended names)} If $d=1,2,3$,
$\sact \in A_1^3$, and the points of the minima belong to the same
straight line, then $(\xb_\ast,t_\ast)$ is a perestroika

$A_1^{+--}$ of the instant shock.

\smallskip\noindent
If $d=2,3$, $\sact \in A_1^4$, and the points of the minima belong
to the same plane, but neither any three of them are on the same
straight line nor all of them are on the same circle, then
$(\xb_\ast,t_\ast)$ is a perestroika

$A_1^{+---}$ provided that the convex hull of our four points is a
triangle;

$A_1^{++--}$ provided that the convex hull of our four points is a
quadrangle.

\smallskip\noindent
If $d=3$, $\sact \in A_1^5$, and neither any four of the points of
the minima belong to the same plane nor all five of them belong to
the same sphere, then $(\xb_\ast,t_\ast)$ is a perestroika $A_1^5$
of the instant shock. More precisely, let us consider the ball of
the minimal radius containing all our five points. Four of them
lie on the boundary of the ball and generate four planes which
divide it into 11 parts: the 1 inner tetrahedron, the 4 domes
adjoining the faces of the tetrahedron, and the 6 lobules
adjoining its edges. We get:

$A_1^{+----}$ if the fifth point is inside of the tetrahedron;

$A_1^{++---}$ if the fifth point is inside of a dome;

$A_1^{+++--}$ if the fifth point is inside of a lobule.

\smallskip\noindent
If $d=1,2,3$, $\sact \in A_3$ such that $\Cb = 0$ when $d=2$ or
$3$, then $(\xb_\ast,t_\ast)$ is a perestroika

$A_3^{-}$ of the instant shock.

\smallskip\noindent
If $d=2,3$, $\sact \in A_1 A_3$, at the minimum $A_3$ $\Cb \ne 0$,
the minimum point $A_1$ belongs to the kernel plane of the minimum
$A_3$ but does not lie either on the kernel line or on the kernel
circle, then $(\xb_\ast, t_\ast)$ is a perestroika $A_1A_3$ of the
instant shock. More precisely, we get (see Figure\,\ref{k}):

$A_1^+A_3^-$ if the point $A_1$ belongs to the disk $\mathcal{D}$;

$A_1^-A_3^+$ if the point $A_1$ belongs to the domain
$\mathcal{U}$;

$A_1^-A_3^-$ if the point $A_1$ belongs to the semiplane
$\mathcal{P}$.

\smallskip\noindent
If $d=3$, $\sact \in A_1^2A_3$, at the minimum $A_3$ $\Cb \ne 0$,
none of the points of the minima $A_1$ lies in the kernel plane of
the minimum $A_3$, the both points $A_1$ and the kernel line do
not belong to the same plane, and the both points $A_1$ and the
kernel circle do not belong to the same sphere, then $(\xb_\ast,
t_\ast)$ is a perestroika $A_1^2A_3$ of the instant shock. More
precisely, let us consider the sphere passing through the both
points $A_1$ and being tangent to the kernel line and the kernel
circle at their common point $O$. This sphere is tangent to the
kernel plane or intersects it along a circle $\mathcal{S}$. When
the kernel plane is between the points $A_1$ we get (see
Figure\,\ref{k}):

$A_1^{++}A_3^-$ if $\mathcal{S}$ lies in the disk $\mathcal{D}$;

$A_1^{--}A_3^+$ if $\mathcal{S}$ lies in the domain $\mathcal{U}$;

$A_1^{--}A_3^-$ if $\mathcal{S}$ lies in the semiplane
$\mathcal{P}$.

\noindent Otherwise, when the points $A_1$ are on one side of the
kernel plane:

$A_1^{+-}A_3^-$ if $\mathcal{S}$ lies in the disk $\mathcal{D}$ or
in the semiplane $\mathcal{P}$, or degenerates into the point $O$
in the case of tangency of the sphere and the kernel plane;

$A_1^{+-}A_3^+$ if $S$ lies in the domain $\mathcal{U}$.

\smallskip\noindent
If $d=3$, $\sact \in A_5$ such that
$$
\left|
\begin{array}{ccc}
\C_1&D_1\\
\C_2&D_2
\end{array}
\right| \ne 0 \quad \mbox{and} \quad \left|
\begin{array}{ccc}
\C_1&D_1&E_1\\
\C_2&D_2&E_2\\
1&0&\C_1^2+\C_2^2
\end{array}
\right| \ne 0,
$$
then $(\xb_\ast,t_\ast)$ is a perestroika

$A_5^{+}$ provided that the indicated determinants have different
signs;

$A_5^{-}$ provided that the indicated determinants have the same
sign.
\end{theorem}
\medskip

As we saw before, the singularities with the same name can occur
at isolated points of the world shocks -- in this case in order to
recognize the perestroika generated by each of them it is enough
to know the extended name of the perestroika -- in this case its
signature is trivial $[0,0]$. Otherwise, the singularities can
form a submanifold $X$ of a positive dimension. In this case the
perestroikas generated by them occur at the points where $X$ is
tangent to isochrones. In order to recognize such a perestroika it
is necessary to know its extended name and signature which is the
signature of the critical point of the restriction $t|_X$ of time
onto $X$. (This critical point is non-degenerate if the initial
condition is generic.) The dimension of $X$ can be calculated
through the multiplicities of the competing global minima:
$$\dim{X}=d+1-\kappa, \quad \kappa = \mu_1+ \dots \mu_\nm -1$$
where $\kappa$ is, of course, just the codimension of $X$.

\medskip
\begin{theorem}
\label{pers} {\sc (Perestroikas -- Signatures)} If $\dim X =0$
then the signature of a perestroika is trivial $[0,0]$. Otherwise,
it is equal to the signature of a quadratic form $\Xi$ on
$\R^{\dim{X}}$ which must not be degenerate in order to the
perestroika to be generic. The form $\Xi$ is defined below.
\end{theorem}

\subsection{Cases $A_1^3$ and $A_1^4$}
\label{qfa1}

Let $d+2 > \nm \ge 3$ and the function $\sact$ have $\nm$ global
minima $A_1$ at points $\Ab_1, \dots, \Ab_\nm$ which lie in the
same $(\nm-2)$-plane $\Lambda$ but do not belong to the same
$(\nm-3)$-dimensional sphere or plane. Then the system of $d+1$
linear equations and one inequality
$$
\left\{
\begin{array}{rcrcr}
           c_1 &+\dots+&              c_\nm &=& 0 \\
 \Ab_1     c_1 &+\dots+&  \Ab_\nm     c_\nm &=& 0 \\
|\Ab_1|^2  c_1 &+\dots+& |\Ab_\nm|^2  c_\nm &<& 0
\end{array}
\right.
$$
has a solution $\mathbf{c}$. Let $\Lambda^\bot$ denote the
$(d+2-\nm)$-dimensional linear space of all (co)vectors being
normal to $\Lambda$ and
$$
\Xi=c_1 \Phi_1 |_{\Lambda^\bot} + \dots + c_m \Phi_m
|_{\Lambda^\bot}
$$
where $\mathbf{c}$ is a solution of our system and $\Phi_i$ is the
quadratic form which is inverse to the matrix of the second
derivatives of $\sact$ at the minimum point $\Ab_i$:
$$
\Phi_i=\Tr_i^{-1}, \quad \Tr_i=
\left\|\sact_{\ab\ab}(\Ab_i)\right\|.
$$

The signature of the quadratic form $\Xi: \Lambda^\bot \to \R$ is
the signature of the perestroika.

\subsection{Case $A_3$}
\label{qfa3}

Let the function $\sact$ have a minimum $A_3$, and $(a,\bb)$ be
orthonormal coordinates giving its representation (\ref{ra3}). We
will need the following extra terms from there:
\begin{equation}
\label{ra3c}
\begin{array}{l}
\sact(a,\bb)=K a^4 +{}\\
\phantom{AAA} + L_1 a^3 \vi_1 + \dots + L_{d-1} a^3 \vi_{d-1} + {}\\
\phantom{AAA}+\Omega_0(\vb) + a \, \Omega_1(\vb) + a^2
\Omega_2(\vb) + \dots
\end{array}
\end{equation}
where $\vb=\bb-\Cb a^2$, $\Omega_0=\omega$, $\Omega_1$, and
$\Omega_2$ are quadratic forms on $\R^{d-1}$ or just symmetric
$(d-1)\times(d-1)$-matrices of real numbers. Let $\Omega =
\Omega_0 + a \, \Omega_1 + a^2 \Omega_2$,
$$
\Tr=
\left\|
\begin{array}{cccc}
4K & 3L_1a^2 & \cdots & 3L_{d-1}a^2\\
L_1&&&\\
\vdots&&2\Omega&\\
L_{d-1}&&&
\end{array}
\right\|
$$
and $\Psi_0$, $\Psi_1$, and $\Psi_2$ be
$(d-1)\times(d-1)$-matrices of real numbers such that
$$
\Tr^{-1} = \left\|
\begin{array}{cc}
\cdot &  \cdots \\
\vdots&\Psi
\end{array}
\right\| + O_{>2}(a)
$$
where  $\Psi = \Psi_0 + a \, \Psi_1 + a^2 \Psi_2$. It is easy to
see that $\Psi$ is symmetric because $\Omega$ is symmetric.

Let $d > 1$ and the function $\sact$ have a global minimum $A_3$
such that $\Cb=0$. The signature of the matrix
$$\Xi=-\Psi_2$$
is the signature of the perestroika.

\subsection{Case $A_1A_3$}

Let $d=3$, $\sact \in A_1A_3$, at the minimum $A_3$ $\Cb \ne 0$,
the minimum point $A_1$ lie in the kernel plane $\Lambda$ of the
minimum $A_3$ but do not belong to the kernel circle. If $\ab =
(a,\bb)$ are orthonormal coordinates giving for $\sact$ the
representation (\ref{ra3c}) and $(a_1, s_1 \, \Cb / |\Cb|)$ are
the coordinates of the minimum point $A_1$ then the system of
linear equations and one inequality
$$
\left\{
\begin{array}{rrcr}
               c_1 &&&{}+c_4  = 0\\
          a_1\,c_1 && {}+c_3& = 0\\
          s_1\,c_1 &{}+|\Cb|\,c_2 && = 0\\
(a_1^2+s_1^2)\,c_1 &{}+\phantom{|\Cb|}\,c_2 && < 0
\end{array}
\right.
$$
has a solution $\mathbf{c}$. Let $(0,\mathbf{n})=(0,\C_2,-\C_1)$
be a normal to $\Lambda$ (co)vector and
$$
\Xi=c_1 \Phi (0,\mathbf{n}) + c_2 \Psi_2 (\mathbf{n}) + c_3 \Psi_1
(\mathbf{n}) + c_4 \Psi_0 (\mathbf{n})
$$
where $\mathbf{c}$ is a solution of our system, $\Phi$ is the
quadratic form defined in \ref{qfa1} for the minimum point $A_1$,
and the matrices of the quadratic forms $\Psi_{0,1,2}$ are defined
in \ref{qfa3}.

The sign of the number $\Xi$ is the signature of the perestroika.

\section{The Minimum Representations}

\label{vsHJ}

\subsection{Unforced case}

We can rewrite our equations (\ref{Burgers}) as one equation for
$\potv$:
$$\potv_t + \frac12 |\nbl{\potv}|^2 = \visc \lpl{\potv}.$$
The Hopf--Cole transformation (which was published by Forsyth back
in 1906 \cite{For06}) $\potv = - 2 \visc \ln \theta$ reduces this
equation to the $d$-dimensional heat equation
$$\potv_t=\visc \lpl{\theta}$$
which can be solved  explicitly if there are no boundaries. We
will get
$$
\potv(\xb,t)=d\,\visc\ln (4 \pi \visc t) -
\phantom{AAAAAAAAAAAAAAA}
$$
$$
{}- 2\visc \ln \int\limits_{\R^d} \exp \left\{
-\frac{\act(\ab,\xb,t)}{2\visc} \right\} \df a_1 \dots \df a_d
$$
where $\act$ is defined in (\ref{W}). This implies our minimum
representation (\ref{W}) for the limit solution as the viscosity
vanishes $\visc \to 0$.

\subsection{Arbitrary potential}

In this case we get $\ham[\potv]= \visc \lpl \potv$ where
$$
\ham[\potv](\xb,t) = \potv_t(\xb,t) + \frac12
|\nbl{\potv(\xb,t)}|^2 + U(\xb, t).
$$
A continuous function $\pot$ is called a {\it viscosity solution}
of the Hamilton--Jacobi equation $\ham[\pot]=0$
if and only if for any point $(\xb_\ast,t_\ast)$:
\begin{itemize}
\item
any smooth function $\phi(\xb,t)$ such that $\pot-\phi$ has a
local minimum at $(\xb_\ast,t_\ast)$ satisfies the inequality
$\ham[\phi](\xb_\ast,t_\ast) \ge 0$, and
\item
any smooth function $\phi(\xb,t)$ such that $\pot-\phi$ has a
local maximum at $(\xb_\ast,t_\ast)$ satisfies the inequality
$\ham[\phi](\xb_\ast,t_\ast) \le 0$.
\end{itemize}

It turns out that, according to \cite{CEL84}, there exists a unique
viscosity solution of the Hamilton--Jacobi equation $\ham[\pot]=0$
with an initial condition $\pot_0$. But the limit solution
$$\pot(\xb,t) = \lim_{\visc \to 0} \potv (\xb,t)$$
satisfies the two above conditions and, therefore, is the
viscosity one. (This fact is a motivation behind the above
definition and term ``viscosity solution''.) Let
$$
F(\ub,\xb,t) = \pot_0(\gamma(0)) + \int\limits_0^{t}
\Lag(\dot{\gamma} (\tau), \gamma (\tau), \tau)\, \df\tau
$$
where $\gamma: [0,t] \to \R^d$ is the solution of the Cauchy
problem:
$$
\left\{
\begin{array}{rl}
\ddot\gamma(\tau) &=-\nbl{U}(\gamma(\tau),\tau)\\
\dot\gamma(t) &= \ub\\
\gamma(t) &= \xb
\end{array}
\right.
$$
and $\Lag(\dot\xb,\xb,t)=|\dot\xb|^2/2-U(\xb,t)$ is the Lagrangian.
According to \cite{ES84}, the function
$$
\pot(\xb,t)=\min_\ub{F(\ub,\xb,t)}
$$
is a unique viscosity solution of the Hamilton--Jacobi
$\ham[\pot]=0$ equation. Therefore it is the minimum representation
for our limit solution.

\subsection{Arbitrary convex Hamiltonian}

Our minimum representation is true for the limit solution of any
equation of the kind:
$$
\potv_t + \Ham(\nbl \potv,\xb, t)=\visc \lpl{\potv}
$$
if the Hamiltonian $\Ham$ is strictly convex down with respect to
momenta. In other words, the matrix of its second derivatives must
be positive definite:
$$
\left\|\Ham_{\pb\pb} (\pb,\xb,t)\right\| > 0.
$$
In the previous case it is true because
$\Ham(\pb,\xb,t)=|\pb|^2/2+U(\xb,t)$.

Let the Lagrangian be the Legendre transformation of the
Hamiltonian:
$$
\Lag(\dot{\xb},\xb,t) = \max_\pb{\left\{\pb \dot{\xb} -
\Ham(\pb,\xb,t) \right\}},
$$
and the trajectory $\gamma$ is the solution of the analogous Cauchy
problem for the Euler--Lagrange equation
$$
\frac{\df}{\df\tau} \frac{\partial}{\partial\dot\gamma}
\Lag(\dot{\gamma},\gamma,\tau) = \frac{\partial}{\partial\gamma}
\Lag(\dot{\gamma},\gamma,\tau).
$$
Then everything works again. Our classification of perestroikas of
shocks is applicable to this most general situation as well.

When the Hamiltonian is not convex only one-dimensional case is
investigated \cite{IK96}, \cite{IK97}. The situation is more
complicated because there is no a minimum representation like
before. Instead a $\min$-$\max$ representation exists \cite{ES84}
but it does not help.

\section{Minimum Functions}

\label{minfuns}

\subsection{Families of functions}

A generic smooth function can have only non-degenerate local
minima $A_1$. All degenerate minima can be killed by an
arbitrarily small perturbation of the original function. In other
words, degenerate minima are not stable with respect to small
perturbations of the original function. For example,  $\ui^4$ --
the simplest degenerate minimum $A_3$ of one variable -- can be
turned into a non-degenerate minima by the perturbation $\ui^4 +
\varepsilon \ui^2$ where $\varepsilon >0$. However, degenerate
minima become stable with respect to small perturbations if we
consider not functions but families of them. For example, the
family $$F(\ui, \lambda) = \ui^4 + \lambda_1 \ui^2 + \lambda_2
\ui$$ of functions of $\ui$ depending on the parameters $\lambda =
(\lambda_1$, $\lambda_2)$ has the minimum $A_3$ if $\lambda =0$.
Any sufficiently small perturbation $F^\ast$ of $F$ will attain a
degenerate minimum along $\ui$ at some point $\ui=\ui_\ast$ for
some parameter values $\lambda = \lambda_\ast$. In order to find
$\ui_\ast$ and $\lambda_\ast$ it is enough to solve the system of
equations $F_\ui^\ast(\ui_\ast,\lambda_\ast) =
F_{\ui\ui}^\ast(\ui_\ast,\lambda_\ast) =
F_{\ui\ui\ui}^\ast(\ui_\ast,\lambda_\ast) = 0$ which has a unique
solution according to the implicit function theorem (as well as
the same system for $F$).

In general, the degenerate minimum $A_{2\nk+1}$ can be attained by
function depending generically on $\n \ge 2\nk$ parameters, and
can be killed by an arbitrarily small perturbation of the original
family for lesser $\n$. The following Theorem\,\ref{fams}
describes locally all generic families containing minima
$A_{2\nk+1}$ of one variable.

\medskip
\begin{lemma}
\label{fams} Let  a function $F(\cdot,\lambda_\ast)$ from a
generic family $F$ of functions of one variable depending on $\n
\ge 2\nk$ parameters have a minimum $A_{2\nk+1}$ at a point
$\ui_\ast$. Then the family $F$ can be reduced in a neighborhood
of the point $(\ui_\ast,\lambda_\ast)$ to the following normal
form:
$$
F(\ui, \lambda) = \ui^{2\nk+2} + \lambda_1 \ui^{2\nk} + \dots +
\lambda_{2\nk} \ui
$$
by a smooth change $\lambda \mapsto h(\lambda)$ of the parameters,
a smooth change $\ui \mapsto g(\ui,\lambda)$ of the variable
depending smoothly on the parameters, and adding $F(\ui,\lambda)
\mapsto F(\ui,\lambda) +  c(\lambda)$ a smooth function of the
parameters.
\end{lemma}
\medskip

Indeed, after the change $\ui \mapsto \ui+\ui_\ast$ we can write
$$
F(\ui,\lambda)=k(\ui,\lambda) \ui^{2\nk+2} + b_0(\lambda)
\ui^{2\nk+1} +{} \qquad
$$
$$
\qquad {}+b_1(\lambda) \ui^{2\nk} + \dots + b_{2\nk}(\lambda) \ui +
c(\lambda)
$$
where $k(0,\lambda_\ast)=1$ and $b_0(\lambda_\ast)=
b_1(\lambda_\ast)= \dots =b_{2\nk}(\lambda_\ast)$. The change $\ui
\mapsto \ui/\sqrt[2\nk+2]{k(\ui,\lambda)}$ gives $k=1$, the change
$\ui \mapsto \ui - b_0(\ui,\lambda)/(2\nk+2)$ kills $b_0$. It
remains to consider $b_1, \dots, b_{2\nk}$ as new parameters in a
neighborhood of $\lambda_\ast$ (taking into account that $F$ is
generic) and kill $c$.

The minimum $X_9$ occurs in generic families of functions only if
the number of parameters $\n \ge 7$. Namely, the $3$ second
derivatives and the $4$ third derivatives must be equal to $0$ --
we get exactly $7$ extra conditions at the minimum point.

\subsection{Singularities}

Let $F(\ub,\lambda)$ be a family of functions of $\ub=(\ui_1,
\dots, \ui_\na)$ depending smoothly on parameters
$\lambda=(\lambda_1, \dots, \lambda_\n)$. Then the function
$$\varphi(\lambda)=\min_\ub F(\ub,\lambda)$$ is called the {\it minimum
function} of the family $F$.

If the family $F$ is smooth and the function
$F(\cdot,\lambda_\ast)$ has the only global minimum which is not
degenerate then the minimum function $\varphi$ is smooth at a
neighborhood of the point $\lambda_\ast$. If the function
$F(\cdot,\lambda_\ast)$ attains its global minimal value at a few
points or its global minimum is degenerate, then the minimum
function $\varphi$ can have singularity at the point
$\lambda_\ast$.

In order to investigate singularities of the minimum functions of
smooth generic families let us note that, while the number of
parameters $\n \le 6$, it is enough to consider only functions of
one variable. Indeed, if $\n \le 6$ and a function
$F(\cdot,\lambda_\ast)$ attains its minimum at a point $\ub_\ast$
we get $\rnk \left \| F_{\ub\ub}(\ub_\ast,\lambda_\ast) \right \|
\ge \na-1$ because the case $\na-2$ requires at least $7$
parameters. Let, for example,
$$
\det \left\| \frac{\partial^2 F}{\partial(\ui_2,\dots,\ui_\na)^2}
\right\| \ne 0
$$
then
$$\varphi(\lambda) =\min_\ui \widetilde{F} (\ui,\lambda), \quad
\widetilde{F}(\ui_1,\lambda)=\min_{\ui_2,\dots,\ui_\na}
F(\ub,\lambda)$$ where $\widetilde{F}$ is a smooth family of smooth
one-variable functions.

If $F$ is a generic family of smooth functions depending smoothly
on $\n \le 4$ parameters $\lambda=(\lambda_1, \dots, \lambda_\n)$
then only the following combinations of minima on the same level
are possible:
\begin{itemize}
\item
$A_1$;
\item
$A_1^2$ if $\n \ge 1$;
\item
$A_1^3$, $A_3$ if $\n \ge 2$;
\item
$A_1^4$, $A_1A_3$ if $\n \ge 3$;
\item
$A_1^5$, $A_1^2A_3$, $A_5$ if $\n = 4$.
\end{itemize}
For example, the combination $A_1^2A_3$ means that our function has two minima
$A_1$, one minimum $A_3$, and the values of the function at all of them are the
same.

Applying Theorem\,\ref{fams} in all these cases to the last
minimum and considering the values at the other minima $A_1$ as
new parameters we get:

\medskip
\begin{theorem}{\rm \cite{Br77}}
\label{Br} If $F$ is a generic family of smooth functions
depending smoothly on $\n \le 4$ parameters
$\lambda=(\lambda_1,\dots,\lambda_\n)$, then its minimum function
$$\varphi(\lambda)=\min_\ui F(\ui,\lambda)$$
can be reduced in a neighborhood of any point $\lambda_\ast$ to
one of the following normal forms such that $\kappa \le \n$:
\begin{itemize}
\item
$A_1^{\kappa+1}$, $\kappa \ge 0$: $\min \{\lambda_1, \dots
\lambda_\kappa,0\}$;
\item
$A_1^{\kappa-2}A_3$, $\kappa \ge 2$:

$\min \left\{\lambda_1, \dots \lambda_{\kappa-2}, \min\limits_\ui
(\ui^4 + \lambda_{\kappa-1} \ui^2 + \lambda_\kappa \ui) \right\}$;
\item
$A_5$, $\kappa=4$: $\min\limits_\ui (\ui^6 + \lambda_1 \ui^4 +
\lambda_2 \ui^3 + \lambda_3 \ui^2 + \lambda_4 \ui)$;
\end{itemize}
by a smooth change of the parameters and adding a smooth function
of them.
\end{theorem}
\medskip

These normal forms are given in Table\,\ref{t} where
$(\alpha,\beta,\gamma,\delta)=(\lambda_1,\dots,\lambda_4)$. All of
them except $A_1$ are not smooth at some points which form a shock
in $\lambda$-space. The shock is smooth in the case
$A_1^2$($\lambda_1=0$), has a boundary in the case $A_3$
($\lambda_1 \le 0, \lambda_2=0$), and has more complicated
singularities in the other cases. Each of the shocks is a cylinder
with $(\n-\kappa)$-dimensional generators. The ``most'' singular
points are defined by the equations
$\lambda_1=\dots=\lambda_\kappa=0$ and are said to be the
corresponding type ($A_1^{\kappa+1}$, $A_1^{\kappa-2} A_3$, or
$A_5$). Any other point has a type with lesser $\kappa$.

These shocks are interpreted as world shocks if $\n=d+1$ and
instant shocks at typical times if $\n=d$ provided that the action
$F$ is a smooth generic function. But instant shocks can
experience perestroikas at separate times which has not been
described by us yet.

\subsection{Perestroikas}
\label{pmf}

Let $F(\ub,\xb,t)$ be a family of functions of $\ub$ depending
smoothly on a point $\xb \in \R^d$ of space and time $t$,
$$
\varphi(\xb,t)= \min_\ub F(\ub,\xb,t),
$$
$$
\shk(\varphi)=\left\{(\xb,t) \in \R^{d+1} \, | \, \varphi \notin
C^\infty(\xb,t) \right\}
$$
be the minimum function and its world shock.

In order to investigate perestroikas of the instant shocks let us
consider the following equivalence of two world shocks:
$$
\shk(\varphi^\prime) \sim \shk(\varphi) \quad  \Longleftrightarrow
\quad \shk(\varphi^\prime) = \g(\shk(\varphi))
$$
where $$\g : \R^{d+1} \to \R^{d+1}, \; (\cdot,t) \mapsto
(\cdot,t+t_\ast)$$ is a sectionally continuously differentiable
homeomorphism of space-time shifting time and being a smooth
diffeomorphism outside the world shocks themselves.

\medskip
\begin{remark}
Unfortunately, if we assume that $\g$ is a smooth diffeomorphism
everywhere we get an infinite local classification of perestroikas
of shocks even if $F$ is generic and $d=2$. Allowing smooth
changes of time instead shifts does not help. For example, the
instant shock at the time of the second perestroika $A_1^4$ in
plane (see Figure\,\ref{f2}) will have at least one continuous
invariant -- the cross-ratio of the four rays.
\end{remark}
\medskip

In order to obtain $\g$ let us consider the graph of the minimum
function $\varphi$
$$
\Gamma(\varphi)=\left\{(\y,\xb,t) \in \R \times \R^{d+1} \, | \,
\y=\varphi(\xb,t) \right\}
$$
and define the equivalence of two graphs:
$$
\Gamma(\varphi^\prime) \sim \Gamma(\varphi) \quad
\Longleftrightarrow \quad \Gamma(\varphi^\prime) =
\G(\Gamma(\varphi))
$$
where
$$
\G : \R \times \R^{d+1} \to \R \times \R^{d+1}, \; (\cdot,\cdot,t)
\mapsto (\cdot,\cdot,t+t_\ast)
$$
is a smooth diffeomorphism shifting time. Then
$$
\g= \pi \circ \G \circ \widetilde{\varphi},
$$
$$
\widetilde{\varphi}: (\xb,t) \mapsto (\varphi(\xb,t),\xb,t), \quad
\pi: (\y,\xb,t) \mapsto (\xb,t).$$

\medskip
\begin{theorem}{\rm \cite{Bog89}}
\label{tb} If $F$ is a generic smooth family and $d \le 3$ then
all perestroikas of instant shocks are locally described by
Table\,\ref{t} with respect to the above equivalence.
\end{theorem}
\medskip

Theorem\,\ref{Br} provides normal forms of world shocks if we set
$\n=d+1$, but the reducing changes of variables mixes space
coordinate and time. In other words, $\lambda$-coordinates from
Theorem\,\ref{Br} depend on both  $\xb$ and $t$. After such mixing
we can represent time as a generic smooth function of $\lambda$
without critical points. In the cases of Theorem\,\ref{Br} such
function $t$ describes a perestroika of the instant shock if its
restriction on the set of the ``most'' singular points
$\lambda_1=\dots=\lambda_\kappa$ has a critical point or
$\kappa=d+1$. Shifting the coordinates
$(\lambda_{\kappa+1},\dots,\lambda_{d+1})$ we get
$$
\begin{array}{c}
t(\lambda)=t_\ast+\tau(\lambda),\\
\tau(0)=0, \quad
\tau_{\lambda_{\kappa+1}}(0)=\dots=\tau_{\lambda_{d+1}}(0)=0.
\end{array}
$$

The graph of the minimum function $\varphi$ is a subset of the
{\it front} $\Sigma(F)$ generated by the family $F$:
$$
\Gamma(\varphi) \subset \Sigma(F) = \bigl\{(\y,\xb,t) \in \R
\times \R^{d+1} \, | \, \exists \, \ub : \phantom{AAAAA}
$$
$$
\phantom{AAAAAAAAA} \y=F(\ub,\xb,t), \, F_\ub(\ub,\xb,t)=0
\bigr\}.
$$
If the function $F(\cdot,\xb_\ast,t_\ast)$ has $\nm$ global minima
with the same minimal value $\y_\ast$ then in a neighborhood of
the point $(\y_\ast,\xb_\ast,t_\ast)$ the front $\Sigma(F)$
consists of $\nm$ branches. In the cases of Theorem\,\ref{Br} we
can, mixing $\y$ and $\lambda$, define $\Sigma(F)$ in the
following more symmetric way:
$$
\begin{array}{l}
\Sigma(F)={}\\
{}=\left\{\left.(\xi,\eta) \in \R^{\kappa+1} \times \R^\nl \right|
\exists \ui : P(\ui)=P^\prime(\ui)=0\right\}
\end{array}
$$
where $\nl=d+1-\kappa$ and in the cases
\begin{itemize}
\item
$A_1^{\kappa+1}$: $\nm=\kappa+1$ and $P=\xi_1 \dots
\xi_{\kappa+1}$;
\item
$A_1^{\kappa-2}A_3$: $\nm=\kappa-1$ and

$P=\xi_1 \dots \xi_{\kappa-2}(\ui^4 + \xi_{\kappa-1} \ui^2
+\xi_\kappa \ui +\xi_{\kappa+1})$;
\item
$A_5$, $\kappa=4$: $\nm=1$ and

$P=\ui^6 + \xi_1 \ui^4 +\xi_2 \ui^3 +\xi_3 \ui^2 +\xi_4 \ui +
\xi_5$.
\end{itemize}
Namely, $\xi_1=\lambda_1-\y, \dots,
\xi_{\nm-1}=\lambda_{\nm-1}-\y$, $\xi_\nm=\lambda_\nm, \dots,
\xi_\kappa=\lambda_\kappa$, $\xi_{\kappa+1}=\y$,
$\eta_1=\lambda_{\kappa+1}, \dots, \eta_\nl=\lambda_{d+1}$.

Let $T: \R^{\kappa+1} \times \R^\nl \to \R$ be a smooth function
without critical points such that
\begin{equation}
\label{T1} T_{\xi_1}(0) \ne 0, \dots, T_{\xi_\nm}(0) \ne 0;
\end{equation}
and if $\nl > 0$ then its restriction onto the plane $\xi=0$ has a
non-degenerate critical point at $0$:
\begin{equation}
\label{T2} T_{\eta_1}(0)=\dots=T_{\eta_\nl}(0)=0, \quad
\det{\left\|T_{\eta\eta}(0)\right\|} \ne 0.
\end{equation}

The {\it extended name} of $T$ codes one of the normal forms
$A_1^{\kappa+1}$, $A_1^{\kappa-2}A_3$, or $A_5$ for $\Sigma(F)$
and the signs of the first derivatives (\ref{T1}). For example:
$A_1^{++---}$ means $A_1^5$ and indicates that among the $5$
numbers (\ref{T1}) there are $2$ positive and $3$ negative ones
(the order is not important -- in extended names pluses are always
followed by minuses); $A_1^{+-}A_3^+$ means $A_1^2A_3$ and
indicates that one of the numbers $T_{\xi_1}(0)$ and
$T_{\xi_2}(0)$ is positive, the other is negative, and the number
$T_{\xi_{3}}(0)$ is positive; $A_5^-$ means $A_5$ and indicates
that $T_{\xi_1}(0) < 0$.

If $\nl \ne 0$ the {\it signature} of $T$ is the signature of the
second differential (\ref{T2}) and is denoted by $[p,q]$ where
$p+q=\nl$.

Here are a few examples. In the case $A_1^2$, $\nl=2$ the function
$-\xi_1 +\xi_2 -\eta_1 \eta_2$ has extended name $A_1^{+-}$ and
signature $[1,1]$; in the case $A_3$, $\nl=1$ the function
$\xi_1-\xi_3 - \eta_1^2$ has extended name $A_3^+$ and signature
$[0,1]$; in the case $A_1A_3$, $\nl=1$ the function $-\xi_1-
\xi_2+\eta_1^2$ has extended name $A_1^-A_3^-$ and signature
$[1,0]$; in the case $A_5$, $\nl=0$ the function has extended name
$A_5^-$.

Theorem\,\ref{tb} follows now from Theorem\,\ref{Arn} proved in
\cite{Arn76}.

\medskip
\begin{theorem}{\rm \cite{Arn76}}
\label{Arn} Let two smooth functions do not have critical points
with the same critical value, satisfy (\ref{T1}) and (\ref{T2}),
and have the same extended name and the same signature. Then they
are translated to each other by a smooth diffeomorphism preserving
the normal form of $\Sigma(F)$.
\end{theorem}
\medskip

Indeed, the translations from Theorem\,\ref{Arn} is performed by
vector fields being tangent to $\Sigma(F)$ (see \cite{Arn76}).
They are, of course, tangent to $\Gamma(F)$ as well. Let $T=\tau
\circ \pi$. The extended name of $T$ can be arbitrary except
$A_1^{+ \dots +}$ and $A_1^{- \dots -}$ because in the case
$A_1^\nm$ $T_{\xi_1} + \dots + T_{\xi_\nm} = 0$ ($\y=\xi_1+ \dots
+ \xi_\nm$ and $T$ does not depend on $\y$). But the normal forms
$\tau$ enumerated in Table\,\ref{t} give for $T$ all possible
combinations of signatures and extended names with these
exceptions.

\subsection{Case where $F$ is an action}

If $F$ is the action generated by a fixed convex smooth Lagrangian
$\Lag$ then the front $\Sigma(F)$ is the many-valued solution of
the Hamilton--Jacobi equation with the initial condition $\pot_0$.
The Hamiltonian is smooth because it is the Legendre
transformation of the Lagrangian -- see Section\,\ref{vsHJ}.
Therefore, if the initial condition $\pot_0$ is generic then the
conditions (\ref{T1}) and (\ref{T2}) are satisfied -- this is
proved in \cite{Bog95} for many-valued solutions of the
Hamilton--Jacobi equation with any smooth Hamiltonian. Hence,
Theorem\,\ref{tb} is true in this situation as well.

\section{The Topological Restrictions}

\label{toprstr}

Let $(\xb_\ast,t_\ast) \in \R^{d+1}$ be an isolated perestroika of
an instant shock. It means the point $(\xb_\ast,t_\ast)$ has a
neighborhood $\mathcal{V} \subset \R^{d+1}$ such that other
perestroikas do not occur in $\mathcal{V}$.

\medskip
\begin{theorem}
\label{myrstr} For a sufficiently small open ball $\mathcal{B}
\subset \R^d$ containing $\xb_\ast$ and for some $\varepsilon > 0$
the local instant shocks $\shk_t \cap \mathcal{B}$ are
contractible (homotopically equivalent to a point) for all $t \in
[t_\ast, t_\ast+\varepsilon)$.
\end{theorem}
\medskip

Originally Theorem\,\ref{myrstr} was discovered and proved in
1989, \cite{Bog89} only for the generic perestroikas when $d \le
3$. But the proof was not topological. Namely, the restriction was
algebraic and forbade some perestroikas from the lists of generic
ones. It turned out that these perestroikas were exactly the same
which were forbidden by the above topological restriction. In
order to explain this fact topologically, Baryshnikov in 1990,
\cite{Bar90}, \cite{ABB91} suggested and proved another
restriction which is true for all isolated (at least) perestroikas
in all dimensions:

\medskip
\begin{theorem}{\rm \cite{Bar90}}
\label{Yurstr} For a sufficiently small open ball $\mathcal{B}
\subset \R^d$ containing $\xb_\ast$ and for some $\varepsilon > 0$
the local complements $\mathcal{B} \setminus \shk_t$ to the
instant shocks are homotopically equivalent to each other for all
$t \in [t_\ast, t_\ast+\varepsilon)$.
\end{theorem}
\medskip

It is easy to check that for the generic perestroikas if $d \le 3$
these two restrictions are equivalent to each other. As a matter
of fact they are equivalent for all isolated perestroikas. Indeed,
roughly speaking the shocks themselves are homotopically
equivalent if and only if their complements are homotopically
equivalent. For example, for the homologies it directly follows
from the Alexander duality between $(\shk_t \cap \mathcal{B}) \cup
\partial{\mathcal{B}}$ and $\mathcal{B} \setminus \shk_t$ where
$\partial{\mathcal{B}}$ is the ball's boundary. But the local
instant shock $\shk_{t_\ast} \cap \mathcal{B}$ at the
perestroika's time is contractible for a sufficiently small ball
$\mathcal{B}\,$!

Below we explain why the homology and homotopy groups of the
complements of $\mathcal{B} \setminus \shk_t$ are the same for $t
\in [t_\ast, t_\ast+\varepsilon)$. These facts follow from
Theorem\,\ref{Yurstr} and are sufficient for all applications. As
a matter of fact, we get weak homotopically equivalences which
imply strong ones for $CW$-complexes. All details can be found in
\cite{Bar90}.

We start with the following simple fact. Let
$$
\mathcal{U}_1 \subset \mathcal{U}_2 \subset \dots, \quad
\mathcal{U}_\infty = \bigcup_{i=1}^\infty \mathcal{U}_i
$$
be an infinite sequence of open subsets in $\R^d$ such that all
embedding $\mathcal{U}_i \hookrightarrow \mathcal{U}_{i+1}$ induce
the isomorphisms
$$
H_\star(\mathcal{U}_i) \cong H_\star(\mathcal{U}_{i+1}), \quad
\pi_\star(\mathcal{U}_i) \cong \pi_\star(\mathcal{U}_{i+1})
$$
between the homology and homotopy groups. But their definitions
use only compact sets and our sets are open, so we get the
isomorphisms
$$
H_\star(\mathcal{U}_i) \cong H_\star(\mathcal{U}_\infty), \quad
\pi_\star(\mathcal{U}_i) \cong \pi_\star(\mathcal{U}_\infty).
$$
If the embeddings do not induce isomorphisms than we can write:
$$
\overrightarrow{\lim} H_\star(\mathcal{U}_i) =
H_\star(\mathcal{U}_\infty), \quad \overrightarrow{\lim}
\pi_\star(\mathcal{U}_i) = \pi_\star(\mathcal{U}_\infty).
$$

Let $\mathcal{U}_t$ denote the set of the Lagrangian coordinates
of all particles which are outside the instant shock $\shk_t$ at a
time $t$. This set is open, diffeomorphic to the complement $\R^d
\setminus \shk_t$ to the shock, and decreases:
$$
\mathcal{U}_{t_1} \subset \mathcal{U}_{t_2} \quad \mbox{if} \quad
t_1 > t_2
$$
because a particle cannot leave the shock. Besides,
$$
\mathcal{U}_{t_\ast} = \bigcup_{t > t_\ast} \mathcal{U}_t.
$$
Now let us assume that there are no perestroikas at all for $t \in
(t_\ast, t_\ast+\varepsilon)$. It means that for the embeddings
$\mathcal{U}_{t_1} \hookrightarrow \mathcal{U}_{t_2}$ induce
isomorphisms
$$
H_\star(\mathcal{U}_{t_1}) \cong H_\star(\mathcal{U}_{t_2}), \quad
\pi_\star(\mathcal{U}_{t_1}) \cong \pi_\star(\mathcal{U}_{t_2})
$$
between the homologies and homotopies for $t_\ast+\varepsilon >
t_1 > t_2 > t_\ast$. Analogously to the above fact for sequences,
we get the isomorphisms:
$$
H_\star(\mathcal{U}_t) \cong H_\star(\mathcal{U}_{t_\ast}), \quad
\pi_\star(\mathcal{U}_t) \cong \pi_\star(\mathcal{U}_{t_\ast})
$$
for $t \in (t_\ast, t_\ast+\varepsilon)$.

For an isolated perestroika we can consider the time $t_\ast$ as
initial -- it means that the Lagrangian coordinates of particles
are defined at $t_\ast$. Besides, we take the sets $\mathcal{U}_t
\cap \mathcal{B}$ instead of $\mathcal{U}_t$.

\section{Shock Recognition}

\label{shkrec}

For the unforced Burgers equation there is the representation
(\ref{W}). Let the function $\sact(\ab)=\act(\ab,\xb_\ast,t_\ast)$
have a collection
$$
\X=\left\{ A_{\mu_1} \dots A_{\mu_\nm} \right\}
$$
of $\nm$ global minima, $\y_\ast$ be the minimal value, and
$\theta(\ab)=\sact(\ab)-\y_\ast$. Let $\E$ be the algebra of
smooth functions on $\R^d$ and
$$
\K(\X) \subset \E
$$
be the ``submanifold'' of all functions having the same collection
$\X$ of global minima such that the values at them are equal to
$0$. In particular, $\theta \in \K(\X)$.

Let us fix the initial condition $\pot_0$ and consider the set of
all functions from the family $\act - \y$:
$$
\W = \left\{ \act(\cdot, \xb,t)-\y \in \E \, | \, \xb \in \R^d, \,
\y, t \in \R, \, t > 0 \right\}.
$$
It is a $(d+2)$-dimensional semiplane:
$$
\W = \pot_0(\ab) + \langle 1, a_1, \dots, a_d \rangle_\R + \langle
|\ab|^2 \rangle_{\R^+} \subset \E
$$
where $\R^+$ is the set of positive real numbers. Moreover, if we
fix $t_\ast$ and consider the corresponding functions from $\act -
\y$:
$$
\W_{t_\ast} = \left\{ \act(\cdot, \xb,t_\ast)-\y \in \E \, | \,
\xb \in \R^d, \, \y \in \R \right\}
$$
then we get a $(d+1)$-dimensional plane
$$
\W_{t_\ast} = \theta(\ab) + \langle 1, a_1, \dots, a_d \rangle_\R
\subset \W.
$$

Let the semiplane $\W$ be transversal to $\K(\X)$ at the point
$\theta$:
$$
\T_\theta\W + \T_\theta\K(\X) = \E
$$
$$
\langle 1, a_1, \dots, a_d, |\ab|^2 \rangle_\R + \T_\theta\K(\X) =
\E.
$$
Then the world shock has the singularity $\X$ at the point
$(\xb_\ast,t_\ast)$. There are two possible cases:

{\it 1) The plane $\W_{t_\ast}$ is transversal to $\K(\X)$ at the
point $\theta$ as well:
$$
\T_\theta\W_{t_\ast} + \T_\theta\K(\X) = \E
$$
$$
\langle 1, a_1, \dots, a_d \rangle_\R + \T_\theta\K(\X) = \E.
$$
Then the instant shock at the time $t_\ast$ has the singularity
$\X$ at the point $\xb_\ast$.

2) The plane $\W_{t_\ast}$ is tangent to $\K(\X)$ at the point
$\theta$. Then the point $(\xb_\ast, t_\ast)$ is a perestroika
$\X$ of the instant shock.}

In order to find its extended name and signature let us note that
the differential of the time coordinate at the point
$(\y_\ast,\xb_\ast,t_\ast)$:
$$
\df t : \T_\theta\W \to \R
$$
defines a unique linear form
$$
\df t : \Q_\theta=\E/\T_\theta\K(\X) \to \R
$$
on the cotangent space to $\K(\X)$ which is uniquely given by its
properties:
\begin{equation}
\label{props}
\begin{array}{c}
\df t([1]) = \df t([a_1])=\dots=\df t([a_d])=0,\\[5pt]
\df t \left(\left[ |\ab|^2 \right]\right) = -2 {t_\ast}^2 < 0
\end{array}
\end{equation}
where $[h] \in \Q_\theta$ denotes the cojugacy class of the
function $h \in \E$. The last property follows from the equation
$\df t \left(\left[ W_t(\ab,\xb_\ast,t_\ast) \right]\right) = 1$.
This linear form is responsible for the extended name of the
perestroika -- it is explained below.

{\it The signature of the perestroika is, by definition, the
signature of the critical point of the restriction
$$
t|_X, \quad X=\W \cap \K(\X);
$$
its second differential in our case is the quadratic form
\begin{equation}
\label{form}
\begin{array}{c}
{\displaystyle \Qd = \df t \circ \II|_{\T_\theta X}, \quad
\T_\theta X = \T_\theta \W_{t_\ast} \cap \T_\theta \K(\X)}\\[5pt]
{\displaystyle \T_\theta \W_{t_\ast} = \langle 1, a_1, \dots, a_d
\rangle_\R}
\end{array}
\end{equation}
where $\II : \T_\theta\K(\X) \to \Q_\theta$ is the second
quadratic form of the ``submanifold'' $\K(\X)$ which is defined
below.}

The ``submanifold'' $\K(\X)$ can be obtained from $\theta$:
$$
\K(\X)=\{\theta(\g(\ab))+h(\ab)\theta(\ab)\}
$$
where $\g:\R^d \to \R^d$ is a smooth diffeomorphism and $h \in
\E$. If $\g_s$ and $h_s$ depend smoothly on a parameter $s \in \R$
we can consider the path
$$
\theta_s(\ab)= \theta(\g_s(\ab)) + h_s(\ab) \theta(\ab) \in
\mathcal{K}_\theta;
$$
we assume that $\g_0(\ab) = \ab$ and $h_0 \equiv 1$, therefore
$\theta_0 = \theta$.

The derivatives of $\g_s$ along $s$ define the vector field
$$
\mathbf{v}_s(\g_s(\ab))=\dot{\g}_s(\ab),
$$
and we can write
\begin{equation}
\label{dtheta}
\dot{\theta}_s(\ab)=(\mathbf{v}_s\cdot\nbl{\theta})(\g_s(\ab)) +
\dot{h}_s(\ab) \theta(\ab).
\end{equation}
In particular,
$$
\dot{\theta}_0=\mathbf{v}_0\cdot\nbl{\theta} + \dot{h}_0 \theta \in
\I_\theta \subset \E
$$
where $\I_\theta$ is the ideal generated by $\theta$ and its all
first partial derivatives. It is the tangent space to $\K(\X)$:
$$\T_\theta \K(\X) = \I_\theta$$
and the factor algebra
$$\Q_\theta=\E/\I_\theta$$
is the cotangent space. Its dimension is equal to the common
multiplicity $\mu=\mu_1+\dots+\mu_\nm$ of the collection $\X$.

Let $\Ab_1, \dots, \Ab_\nm \in \R^d$ be the minimum points
$A_{\mu_1}, \dots, A_{\mu_\nm}$ respectively and
$$
\I_{\theta,i} = \I_\theta + \goth{s}_i, \quad \Q_{\theta,i} =
\E/\I_{\theta,i}
$$
where $\goth{s}_i \subset \E$ is the ideal consisting of all
functions such that each of them is equal to $0$ in some
neighborhood (depending on the function) of the point $\Ab_i$. Then
we get the decompositions
\begin{equation}
\label{dc}
\begin{array}{rcrcccr}
\I_\theta & = & \I_{\theta,1} & \cap & \dots & \cap & \I_{\theta,\nm},\\
\Q_\theta & = & \Q_{\theta,1} & \oplus & \dots & \oplus &
\Q_{\theta,\nm}.
\end{array}
\end{equation}
In each algebra $\Q_{\theta,i}$ there is a unique maximal ideal
$\goth{m}_{\theta,i}$. Its degree $\goth{m}_{\theta,i}^{\mu_i-1}$
defines the so-called Jacobian ideal which is a unique minimal
ideal and a one-dimensional $\R$-linear subspace of the algebra
$\Q_{\theta,i}$. This ideal is generated by the element
$$
\goth{j}_{\theta,i} = \left[ |\ab-\Ab_i|^{\mu_i-1} \right] \in
\Q_{\theta,i} \subset \Q_\theta
$$
($\mu_i$ is odd) and the ray $\langle \goth{j}_{\theta,i}
\rangle_{\R^+} \subset \Q_\theta$ does not depend on coordinates
($\R^+$ is positive numbers).

{\it The signs of the numbers
$$
c_1=\df t(\goth{j}_{\theta,1}), \dots, c_\nm=\df
t(\goth{j}_{\theta,\nm})
$$
form the extended name of perestroika.}

If $\X=A_1^\nm$ then $\mu=\nm$,
$$
\I_{\theta, i} = (a_1-\Ai_{i1}, \dots, a_d-\Ai_{id})
$$
is a maximal ideal where $(\Ai_{i1}, \dots, \Ai_{id})$ are the
coordinates of the minimum point $\Ab_i$,
$$
\Q_{\theta,i} \cong \R, \quad \goth{j}_{\theta,i}=[1], \quad
\Q_\theta \cong \R^\nm,
$$
and the decomposition of $[1], [a_1], \dots, [a_d], [|\ab|^2] \in
\Q_\theta$ can be written as a $(d+2)\times \mu$-matrix:
$$
M = \left\|
\begin{array}{c}
{[1]} \\
{[\ab]}\\
\left[ |\ab|^2 \right]
\end{array}
\right\|=\left\|
\begin{array}{ccc}
1 & \dots & 1 \\
\Ab_1 & \dots & \Ab_\nm\\
|\Ab_1|^2 & \dots & |\Ab_\nm|^2
\end{array}
\right\|.
$$

Let $\X=A_1^lA_3$, $l=\nm-1$ and $(a,\bb)$ be orthonormal
coordinates giving the representation (\ref{ra3}) for the minimum
$A_3$ of $\theta$. Then $\mu=\nm+2$ and using the coordinates
$(a,\vb)$ we get
$$
\I_{\theta,\nm} = (a^3+O_{>3},\vi_1+O_{>2},\dots,\vi_{d-1}+O_{>2});
$$
the Nakayama lemma implies that
$$
\I_{\theta,\nm} = (a^3,\vi_1,\dots,\vi_{d-1}).
$$
Therefore, $\Q_{\theta,\nm} \cong \R[a]/(a^3)$,
$\goth{j}_{\theta,\nm}=[a^2]$,
$$
\quad \Q_\theta \cong \R^{\nm-1} \oplus \R[a]/(a^3).
$$
Using the basis $[a^2],[a],[1]$ in $\R[a]/(a^3)$ the decomposition
of $[1], [a], [b_1], \dots, [b_{d-1}], [|\ab|^2] \in \Q_\theta$ can
be written as a $(d+2)\times \mu$-matrix
$$
M=\left\|
\begin{array}{cccccc}
1 & \dots & 1 & 0 & 0 & 1 \\
\Ai_1 & \dots & \Ai_l & 0 & 1 & 0\\
\Bb_1 & \dots & \Bb_l & \Cb & 0 & 0\\
\Ai_1^2+|\Bb_1|^2 & \dots & \Ai_l^2+|\Bb_l|^2 & 1 & 0 & 0
\end{array}
\right\|
$$
where $(\Ai_i,\Bb_i)$ are the orthonormal coordinates of the
minimum point $\Ab_i$ and $\Cb$ comes from the representation
(\ref{ra3}): $[\bb]=\Cb [a^2] \pmod{\I_{\theta,\nm}}$.

Let $\X=A_5$ and $(a,\bb)$ be orthonormal coordinates giving the
representation (\ref{ra5}) for the minimum $A_5$ of $\theta$. Then
$\nm=1$, $\mu=5$ and using the coordinates $(a,\vb)$ we get
$$
\I_\theta = (a^5+O_{>5},\vi_1+O_{>4},\dots,\vi_{d-1}+O_{>4})
$$
and the Nakayama lemma implies that
$$
\I_\theta = (a^5,\vi_1,\dots,\vi_{d-1}).
$$
Therefore, $\Q_\theta \cong \R[a]/(a^5)$, $\goth{j}_\theta=[a^4]$.

Using the basis $[a^4],\dots,[a],[1]$ in $\R[a]/(a^5)$ the
decomposition of $[1], [a], [b_1], \dots, [b_{d-1}], [|\ab|^2] \in
\Q_\theta$ can be written as a $(d+2)\times \mu$-matrix
$$
M=\left\|
\begin{array}{ccccc}
0          & 0          & 0   & 0 & 1\\
0          & 0          & 0   & 1 & 0\\
\mathbf{E} & \mathbf{D} & \Cb & 0 & 0\\
|\Cb|^2    & 0          & 1   & 0 & 0
\end{array}
\right\|
$$
where $\Cb, \mathbf{D}, \mathbf{E}$ comes from the representation
(\ref{ra5}): $[\bb]=\Cb [a^2] + \mathbf{D} [a^3] + \mathbf{E} [a^4]
\pmod{\I_\theta}$.

Let, in these all cases,
$$
\df t = (c_1, \dots, c_\mu)
$$
be the coordinates of the linear form $\df t$ written as a
$\mu$-column $\mathbf{c}$ and the $(d+1) \times \mu$-matrix
$M^\prime$ be the matrix $M$ without the last row $M_{d+2}$:
$$
M=\left\|
\begin{array}{c}
M^\prime \\ M_{d+2}
\end{array}
\right\|.
$$

So, if $\rnk{M}=\mu$ then the semiplane $\W$ is transversal to
$\K(\X)$ at the point $\theta$ and the world shock has the
singularity $\X$ at the point $(\xb_\ast,t_\ast)$. There are two
possible cases:

{\it 1) If $\rnk{M^\prime}=\mu$ too then the plane $\W_{t_\ast}$
is transversal to $\K(\X)$ at the point $\theta$ as well and the
instant shock at the time $t_\ast$ has the singularity $\X$ at the
point $\xb_\ast$. This is an algebraic form of
Theorem\,\ref{sing}.

2) If $\rnk{M^\prime}=\mu-1$ then the plane $\W_{t_\ast}$ is
tangent to $\K(\X)$ at the point $\theta$ and the point
$(\xb_\ast, t_\ast)$ is a perestroika $\X$ of the instant shock.
Besides, the system of $d+1$ linear equations and one inequality
coming from (\ref{props}):
$$
\left\{
\begin{array}{rcc}
M^\prime \mathbf{c} & = & \mathbf{0} \\
M_{d+2} \mathbf{c} & < & 0
\end{array}
\right.
$$
has a solution and the signs of its components $c_1, \dots, c_\nm$
give the perestroikas's extended name. Neither of them must equal
to $0$ -- otherwise the perestroika is not generic and does not
belong to our lists. This is an algebraic form of
Theorem\,\ref{pern}.}

In order to calculate the perestroika's signature let us define
the second quadratic form of the ``submanifold'' $\K(\X) \subset
\E$ at the point $\theta$:
$$
\II_\theta: \I_\theta \times \I_\theta \to \Q_\theta,
$$
$$
\II_\theta(\goth{x},\goth{y}) = [ \mathbf{u} \cdot \nbl{\goth{y}} ]
= [ \mathbf{v} \cdot \nbl{\goth{x}} ] \pmod{\I_\theta},
$$
$$
\goth{x}=\mathbf{u}\cdot\nbl{\theta} + h \theta \in \I_\theta,
$$
$$
\goth{y}=\mathbf{v}\cdot\nbl{\theta} + k \theta \in \I_\theta.
$$
These formulas show that $\II$ is correctly defined (depends only
on $\goth{x,y}$ but not on $\mathbf{u,v},k,h$), symmetric, and
linear over $\E$. Besides,
$$
\J_\theta=\I_\theta^2+(\theta) \subset \ker{\II_\theta}.
$$
As usual, the geometrical sense of the second quadratic form $\II$
is the curvature of $\K(\X)$:
$$
\II_\theta(\dot{\theta}_0,\dot{\theta}_0)=\ddot{\theta}_0
\pmod{\I_\theta}
$$
because differentiating (\ref{dtheta}) we get
$$
\ddot{\theta}_0 = \mathbf{v}_0 \cdot \nbl {(\mathbf{v}_0 \cdot
\nbl{\theta})} + \dot{\mathbf{v}}_0\cdot\nbl{\theta} + \ddot{h}_0
\theta.
$$
We can decompose the second quadratic form using the decomposition
(\ref{dc}):
$$
\II_\theta = \II_{\theta,1} \oplus \dots \oplus \II_{\theta,\nm},
\quad \II_{\theta,i}: \I_{\theta,i} \to \Q_{\theta,i},
$$
$$
\J_{\theta,i}=\I_{\theta}^2+(\theta)+\goth{s}_i \subset
\ker{\II_{\theta,i}}.
$$

Let $\X=A_1^\nm$, then $\J_{\theta,i}=\I_{\theta,i}^2$. How to
calculate the matrix ${\left(\Phi_i\right)}_{qr}=\II_{\theta,i}
(a_q-\Ai_{iq},a_r-\Ai_{ir})$? We know that
$$
\partial_{a_p} \theta = \sum_q {\left(\Tr_i\right)}_{pq}
(a_q-\Ai_{iq}) \pmod{\J_{\theta,i}},
$$
then the equality
$$
\II_{\theta,i} (\partial_{a_p} \theta, a_r-\Ai_{ir}) =
\partial_{a_p} a_r = \delta_{pr}
$$
gives $ \Phi_i=\Tr_i^{-1} \pmod{\I_{\theta,i}}$. But, according to
(\ref{form}),
$$
\T_\theta X = \langle 1, a_1, \dots, a_d \rangle_\R \cap \I_\theta
= \Lambda^\bot
$$
and $\Xi= c_1 \Phi_1|_{\Lambda^\bot} + \dots + c_m
\Phi_m|_{\Lambda^\bot}$.

Let $\X=A_3$, then $\I_\theta = (a^3,\vi_1,\dots,\vi_{d-1})$,
$\J_\theta=\I_\theta^2 +(a^4)$. How to calculate the matrix
$\Psi^\prime_{qr}=\II_\theta (\goth{x}_q,\goth{x}_r)$ where
$(\goth{x}_1,\dots,\goth{x}_d) = (a^3,\vi_1,\dots,\vi_{d-1})$?
Using the coordinates $\ub=(a,\vb)$ we get
$$
\partial_{\ui_p} \theta =\sum_q \Tr_{pq} \goth{x}_q \pmod
{\J_\theta},
$$
then the equality $\II (\partial_{\ui_p} \theta, \goth{x}_r) =
\partial_{\ui_p} \goth{x}_r$ gives
$$
\Psi^\prime=\Tr^{-1} \left\|
\begin{array}{cccc}
3a^2      & 0\\
0& \id\\
\end{array}
\right\| \pmod{\I_\theta}
$$
where $\id$ is the identity $(d-1)\times(d-1)$-matrix and the last
multiplier is diagonal.

If $\Cb=0$ then $\bb=\vb$ -- see (\ref{ra3}) -- and
$$
\T_\theta X = \langle 1, a, b_1 \dots, b_{d-1} \rangle_\R \cap
\I_\theta = \langle \vi_1 \dots, \vi_{d-1} \rangle_\R.
$$
But $\Psi=\Psi^\prime|_{\langle \vi_1 \dots, \vi_{d-1}
\rangle_\R}$ and, according to (\ref{form}),
$$
\Xi= c_1 \Psi_2 + c_2 \Psi_1 + c_3 \Psi_0 = -\Psi_2
$$
because $\mathbf{c}=(-1,0,0)$.

If $\X=A_1A_3$, $d=3$, $\Cb \ne 0$, and $\Ab_1$ belongs to the
kernel plane $\Lambda$. Then
$$
\T_\theta X = \langle 1, a, b_1, b_2 \rangle_\R \cap \I_\theta =
\langle \C_2 \vi_1 - C_1 \vi_2 \rangle_\R = \Lambda^\bot
$$
and, according to (\ref{form}),
$$
\Xi=c_1 \Phi (0,\mathbf{n}) + c_2 \Psi_2 (\mathbf{n}) + c_3 \Psi_1
(\mathbf{n}) + c_4 \Psi_0 (\mathbf{n})
$$
where $\mathbf{n}=(\C_2,-\C_1)$.

\section*{Acknowledgements}

This work was started during a visit to the Laboratoire
G.\,D.\,Cassini of the Observatoire de la C\^ote d'Azur, supported
by the French Ministry of Education. During this visit I had very
useful discussions with J.\,Bec, M.\,Blank, U.\,Frisch,
G.\,Molchan, and T.\,Matsumoto. The work was completed with
support of the Levelhulme Trust under hospitality of the
Department of Mathematical Sciences of the University of
Liverpool.

\bibliography{shock}

\def\cprime{$'$}
\begin{thebibliography}{10}

\bibitem{FB00}
U.~Frish and J.~Bec.
\newblock Burgulence.
\newblock In M.~Lesieur, editor, {\em Les Housches 2000: New Trends in
  Turbulence}. Springer EDP-Sciences, to appear.
\newblock arXiv:nlin.CD/0012033.

\bibitem{Zel70}
Ya.~B. Zel{\cprime}dovich.
\newblock Gravitational instability: an approximate theory for large density
  perturbations.
\newblock {\em Astron. Astrophys.}, 5:84--89, 1970.

\bibitem{AShZ}
V.~I. Arnol{\cprime}d, S.~F. Shandarin, and Ya.~B. Zel{\cprime}dovich.
\newblock The large-scale structure of the universe. {I}. {G}eneral properties.
  {O}ne and two-dimensional models.
\newblock {\em Geophys. Astrophys. Fluid Dynam.}, 20:111--130, 1982.

\bibitem{ArnCat}
V.~I. Arnol{\cprime}d.
\newblock {\em Catastrophe theory}.
\newblock Springer-Verlag, Berlin, third edition, 1992.
\newblock Translated from the Russian by G. S. Wassermann, Based on a
  translation by R. K. Thomas.

\bibitem{AGVI}
V.~I. Arnol{\cprime}d, S.~M. Guse{\u\i}n-Zade, and A.~N. Varchenko.
\newblock {\em Singularities of differentiable maps. {V}ol. {I}}.
\newblock Birkh\"auser Boston Inc., Boston, MA, 1985.
\newblock The classification of critical points, caustics and wave fronts,
  Translated from the Russian by Ian Porteous and Mark Reynolds.

\bibitem{GMS91}
S.~N. Gurbatov, A.~N. Malakhov, and A.~I. Saichev.
\newblock {\em Nonlinear random waves and turbulence in nondispersive media:
  waves, rays, particles}.
\newblock Manchester University Press, Manchester, 1991.
\newblock Translated from the Russian. Translation edited and with a preface by
  D. G. Crighton.

\bibitem{CEL84}
M.~G. Crandall, L.~C. Evans, and P.-L. Lions.
\newblock Some properties of viscosity solutions of {H}amilton-{J}acobi
  equations.
\newblock {\em Trans. Amer. Math. Soc.}, 282(2):487--502, 1984.

\bibitem{For06}
A.~R. Forsyth.
\newblock {\em Theory of differential equations. Partial differential
  equations}, volume~6, chapter XII, pages 100--102.
\newblock Cambridge University Press, 1906.

\bibitem{Kr66}
S.~N. Kru{\v{z}}kov.
\newblock Generalized solutions of nonlinear equations of the first order with
  several variables. {I}.
\newblock {\em Mat. Sb. (N.S.)}, 70 (112):394--415, 1966.

\bibitem{Kr67}
S.~N. Kru{\v{z}}kov.
\newblock Generalized solutions of nonlinear equations of the first order with
  several independent variables. {I}{I}.
\newblock {\em Mat. Sb. (N.S.)}, 72 (114):108--134, 1967.

\bibitem{Bog99}
I.~A. Bogaevsky.
\newblock Singularities of viscosity solutions of {H}amilton-{J}acobi
  equations.
\newblock {\em S\=urikaisekikenky\=usho K\=oky\=uroku}, (1111):138--143, 1999.
\newblock Singularity theory and differential equations (Kyoto, 1999).

\bibitem{ES84}
L.~C. Evans and P.~E. Souganidis.
\newblock Differential games and representation formulas for solutions of
  {H}amilton-{J}acobi-{I}saacs equations.
\newblock {\em Indiana Univ. Math. J.}, 33(5):773--797, 1984.

\bibitem{Br77}
L.~N. Bryzgalova.
\newblock Singularities of the maximum of a function that depends on the
  parameters.
\newblock {\em Funkcional. Anal. i Prilo\v zen.}, 11(1):59--60, 1977.

\bibitem{Bog89}
I.~A. Bogaevsky.
\newblock Metamorphoses of singularities of minimum functions, and bifurcations
  of shock waves of the {B}urgers equation with vanishing viscosity.
\newblock {\em Algebra i Analiz}, 1(4):1--16, 1989.
\newblock English translation in St.\,Petersburg (Leningrad) Math. J., 1
  (1990), no.\,4, 807--823.

\bibitem{Bog95}
I.~A. Bogaevsky.
\newblock Perestroikas of fronts in evolutionary families.
\newblock {\em Trudy Mat. Inst. Steklov.}, 209:65--83, 1995.
\newblock English translation in Proc. Steklov Inst. Math. 209 (1995), 57--72.

\bibitem{Bar90}
Yu.~M. Baryshnikov.
\newblock Topology of surgeries of nonsmoothness sets of minimum functions of
  variational problems.
\newblock {\em Funktsional. Anal. i Prilozhen.}, 24(3):62--63, 1990.
\newblock English translation in Funct. Anal. Appl. 24 (1990), no.\,3,
  220--221.

\bibitem{IK96}
Sh. Izumiya and G.~T. Kossioris.
\newblock Formation of singularities for viscosity solutions of
  {H}amilton-{J}acobi equations.
\newblock In {\em Singularities and differential equations (Warsaw, 1993)},
  pages 127--148. Polish Acad. Sci., Warsaw, 1996.

\bibitem{IK97}
Sh. Izumiya and G.~T. Kossioris.
\newblock Bifurcations of shock waves for viscosity solutions of
  {H}amilton-{J}acobi equations of one space variable.
\newblock {\em Bull. Sci. Math.}, 121(8):619--667, 1997.

\bibitem{Arn76}
V.~I. Arnol{\cprime}d.
\newblock Wave front evolution and equivariant {M}orse lemma.
\newblock {\em Comm. Pure Appl. Math.}, 29(6):557--582, 1976.

\bibitem{ABB91}
V.~I. Arnol{\cprime}d, Yu.~M. Baryshnikov, and I.~A. Bogaevsky.
\newblock Singularities and bifurcations of potential flows.
\newblock In {\em {\em S.\,N.\,Gurbatov, A.\,N.\,Malakhov, and A.\,I.\,Saichev}
  Nonlinear random waves and turbulence in nondispersive media: waves, rays,
  particles}, pages 290--300. Manchester University Press, Manchester, 1991.

\end{thebibliography}
\bibliographystyle{unsrt}

\end{document}